\documentclass{article}

\usepackage{dirtytalk, enumitem}

\usepackage{amsmath, amsthm, appendix}

\usepackage{mathtools}

\usepackage[section]{algorithm}
\usepackage{algpseudocode}

\usepackage{pgf, tikz}
\usetikzlibrary{angles, quotes, automata, positioning}
\usetikzlibrary{arrows}

\usepackage{subcaption}

\usepackage[margin=1.0in]{geometry}
\usepackage{graphicx}

\usepackage{hyperref}

\theoremstyle{definition}

\algblockdefx[Match]{Match}{EndMatch}
    [1]{\textbf{match} #1}
    {\textbf{end match}}

\algloopdefx[Case]{Case}
    [1]{\textbf{case} #1:}

\algblockdefx[ForIn]{ForIn}{EndForIn}
    [2]{\textbf{for} #1 \textbf{in} #2}
    {\textbf{end for}}

\begin{document}

\title{One Hundred and Twelve Point Three Degree Theorem}

\author{George Tokarsky, Jacob Garber, Boyan Marinov, Kenneth Moore\\University of Alberta}

\maketitle

\section{Introduction}

It has been known since Fagnano in 1775 that an acute triangle always has a periodic billiard path, namely the orthic triangle. In 1993, Holt \cite{holt} showed that every right triangle has a periodic path using simple arguments. It is currently unknown whether every obtuse triangle has a periodic path, and there don't appear to be any simple arguments. In between, there have been various partial results. In 1986, Masur \cite{masur} proved that every rational obtuse triangle has a periodic path, and in 2006, Schwartz \cite{schwartz} showed that every obtuse triangle with obtuse angle at most 100 degrees has a periodic path using a computer assisted proof. The aim of this paper is to show that every obtuse triangle with obtuse angle at most 112.3 degrees has a periodic path using a different computer assisted proof.

\section{Side, Code and Alphabet Sequences}

Drawing inspiration from the game of billiards, consider a frictionless particle moving about the interior of a triangle. If the particle encounters a side, it bounces off according to the law of reflection, and if the particle hits a vertex, it is considered to end there. A \textbf{billiard trajectory} or \textbf{poolshot} is the path this particle traces out as it moves within the triangle. A poolshot is said to be \textbf{periodic} if the particle eventually returns to its starting point and repeats the same path over again.

\subsection{Side Sequences}

To begin classifying billiard trajectories, we first need some notation. Given a triangle ABC with angles \( x \), \( y \), and \( z \), we label the side AB opposite \( z \) with 1, the side BC opposite \( x \) with 2, and the side AC opposite \( y \) with 3. For any periodic billiard trajectory in the triangle, we then define its \textbf{side sequence} to be the list of consecutive sides that are hit during one period. For example, for the periodic billiard trajectory shown in Figure 1, if we begin at side 2 and continue to side 1, it has the side sequence 213132313 of \textbf{length} 9. Observe that no two consecutive positive integers from 1,2,3 are the same including the first and the last which we formally call a \textbf{legal side sequence}. We call a finite side sequence \textbf{repeating} if the last integer in it is followed by the first integer and then the integers keep repeating.

\begin{figure}[ht]
    \centering
    \includegraphics[scale=0.12]{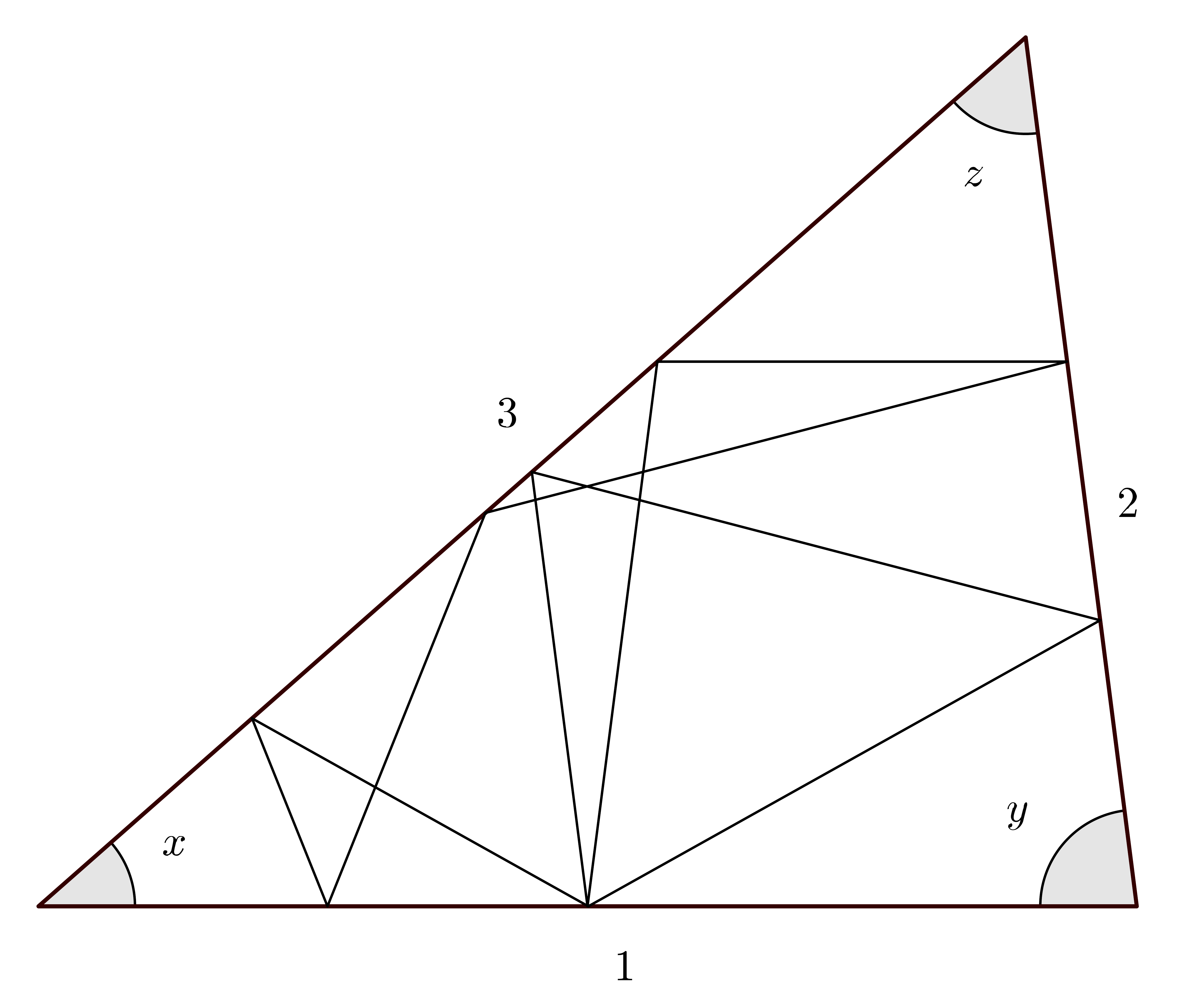}
    \caption{Periodic Billiard Trajectory 213132313}
    \label{fig:traj}
\end{figure}

There are other equivalent side sequences for the trajectory in Figure \ref{fig:traj} that we could consider. Given the periodic nature of the trajectory, the initial starting side and direction are both arbitrary, so we could instead start on side 3 and continue to side 1, giving the side sequence 313213132. These choices correspond to rotations and reversals of the original sequence, respectively. Likewise, we could continue around the triangle for as many periods as we wish, leading to the infinite family of side sequences 313213132313213132, 313213132313213132313213132, etc. All of these periodic sequences represent the same trajectory, so it is convenient to pick one among them as being canonical.

A periodic side sequence is said to be in \textbf{standard form} if it has one period, and is lexicographically least among all its rotations and reversals. All side sequences can be converted to an equivalent one in standard form. For example, the standard form of the side sequence 213132313 is 123132313. In this way, every periodic billiard trajectory in triangle ABC can be uniquely identified by a side sequence in standard form.

A periodic side sequence is said to be in \textbf{extra standard form} if we renumber the sides to make the side sequence minimal. For example the side sequence 3132 can be minimalized to  1213 by renumbered the sides. In this way any side sequence can be written in the form 12...23 or 12...13  

Note that given a legal side sequence, there may not exist any triangle with a periodic billiard trajectory of that type. We call such side sequences \textbf{empty}.
Observe also that any periodic path in a triangle must hit all three sides since it will eventually bounce out of any angle of the triangle and hence any periodic side sequence must include all three integers 1, 2 and 3.

\subsection{Code Sequences}

Long side sequences quickly become cumbersome to work with, so it is helpful to introduce a more compact notation. For each pair of adjacent integers in the side sequence (including the first and last), write out the angle between those corresponding sides. For example in the side sequence 123231323, this is \(y z z z x x z z x \) where the first y is between the first 1 and 2 and right up to the last x coming between the last 3 and the first 1. Then, group the consecutive angles together, shuffling identical angles from the back to the front if necessary, and count the number of angles in each group. This gives \( 1y\ 3z\ 2x\ 2z\ 1x \). The sequence of integers 1 3 2 2 1 with spaces is called the \textbf{code sequence} consisting of the five individual positive integers each called a \textbf{code number}, and the sequence of angles \( y\ z\ x\ z\ x \) with spaces is called the \textbf{angle sequence}. The \textbf{length} of the code is the number of code numbers called an \textbf{odd code} or \textbf{even code} depending on its length and the \textbf{sum} of the code is the sum of the code numbers which is the same as the length of its corresponding side sequence. Given a code and angle sequence, it is possible to recover the original side sequence if we start it with the first two integers 12... as in our example.

Alternately observe the first code number 1 represents how many times the integers 1 and 2 are interchanged in a row at the start of the side sequence, the second code number 3 then represents the following number of succesive interchanges of 2 and 3, the third code number 2 then represents the following number of succesive interchanges of 3 and 1 and so forth. In this process assuming we are dealing with a side sequence of a periodic path, we always view the last integer as followed by the first integer and thus in this example there is just one interchange of 3 and 1 at the end. Also note that any code number by itself plus one yields the number of two successive  alternate integers in the side sequence. For example the code number 2 corresponds to a side sequence of the form aba which is symmetric about the middle b.

Note that specifying any two consecutive angles in the angle sequence will completely determine the rest of it as long as we know the code sequence. For example, for the sequence 1 3 2 2 1, consider the angles \( y\ z \) for the 1 and the 3. After bouncing once across angle \( y \), the pool ball will bounce 3 times across angle \( z \), and 3 being odd, will end at the opposite side from where it started. It will then bounce twice across angle \( x \), and 2 being even, will end up at the same side as where it started, and then bounce back across the previous angle \( z\). Continuing this pattern will produce the rest of the angles in the angle sequence and, crucially, wrap-around to match the last \( x \) with the first \( y \) we started with.

    As with side sequences, a code sequence is in \textbf{standard form} if it is lexicographically least among its rotations and reversals. For example, the standard form of 1 1 3 2 2 is 1 1 2 2 3.

We will call a finite string of positive integers separated by spaces a \textbf{legal code sequence} if it is the \textbf{code sequence of a legal side sequence} and hence potentially represents a periodic path in a triangle.

\textbf{Important comment one:} We will only use the code sequence notation when dealing with a \textbf{repeating} side sequence corresponding to that code sequence. 

\noindent

\textbf{Important comment two:} Observe that if a side sequence starts 13 and ends in 2, then given its corresponding code sequence, we can recover the original side sequence by starting the sequence with 13... since successive integers are completely determined by the code numbers. If we choose to recover the side sequence by starting it with 12... then we will recover it with the 3's and 2's interchanged which means it will now end in a 3. This is not a problem as it just represents a relabeling of the sides of triangle ABC. Similarly we can start the side sequence 23 or any combination of 1,2 and 3 with the corresponding relabelling of the triangle. Caution: If the original side sequence is not in standard form and for example starts 13 and ends in 3, then we will recover some relabelling of the original side sequence from the code.

\textbf{Important comment three:} Given a code sequence representing a legal side sequence then we can always write in the form ab...bc or ab...ac where a,b and c are distinct.

We also make the convention that if there are three dots in front of (or following) a sequence of code numbers then this means there is at least one \textbf{code number} preceding (or following) that sequence in which case we will call it a \textbf{subcode}. For example ...2 4... is a subcode of the code sequence 2 2 4 4. Caution: A subcode need not be a legal code sequence by itself.

\subsection{Alphabet Code Sequences}

We write each code sequence in terms as odd and even integers without spaces.  For instance, the \textbf{alphabet code sequence} of 1 1 2 2 3 is OOEEO. A alphabet code sequence is said to be in \textbf{standard form} if it is alphabetically least among all its rotations and reversals. For example, the standard form of OOEEO is EEOOO.

\noindent
It is very easy to determine if a code sequence is legal or not as in the automaton below.

\noindent
\textbf{Algorithm:} An alphabet code sequence is legal if and only if it forms a closed path of E's and O's as in the automaton.

\noindent
\begin{proof}

Each circle corresponds to a side of a triangle and two angles adjacent to it, say x y and then make a sequence of O's and E's. Every E ends up at the same side but in the opposite direction while every O ends up at the next side in the same direction. This means a periodic path can only start at some side (corresponding circle) and continue in the same direction. Equivalently this means there must be a closed path of O's and E's in the automaton. One can then check that a closed path can then be reduced to the empty set $\phi$ by successively eliminating strings of type EE, OOO, OEOE or EOEO, OOEOOE, OEOOEO or EOOEOO.
\end{proof}

\begin{figure}[ht]
    \centering

\begin{tikzpicture}[->, initial text=, shorten >=1pt, node distance=2cm, auto, semithick, >=stealth]

    \node[state, initial, accepting] (xy) { x y };
    \node[state] (yz) [right=of xy] {y z};
    \node[state] (zx) [right=of yz] {z x};

    \node[state] (yx) [below=of xy] {y x};
    \node[state] (zy) [below=of yz] {z y};
    \node[state] (xz) [below=of zx] {x z};

    \path (xy) edge              node {O} (yz)
               edge [bend left]  node {E} (yx)
          (yz) edge              node {O} (zx)
               edge [bend left]  node {E} (zy)
          (zx) edge [bend right] node [swap] {O} (xy)
               edge [bend left]  node {E} (xz)
          (yx) edge [bend right] node [swap] {O} (xz)
               edge [bend left]  node {E} (xy)
          (zy) edge              node [swap] {O} (yx)
               edge [bend left]  node {E} (yz)
          (xz) edge              node [swap] {O} (zy)
               edge [bend left]  node {E} (zx);
\end{tikzpicture}

    \caption{Automaton for Alphabet Code Sequences}
    \label{fig:dfa}
\end{figure}
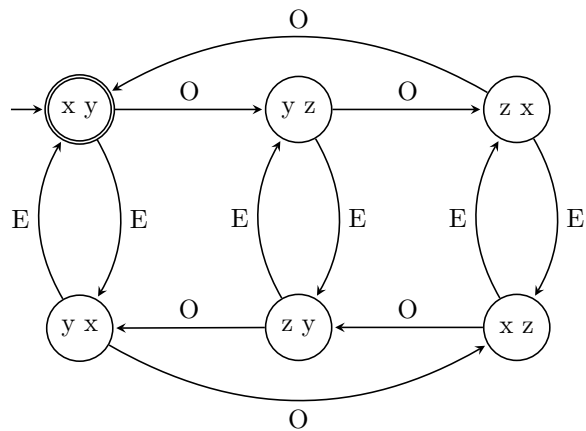

    Observe from the automaton that a code sequence always has an even number of even numbers. It follows that the length and sum of a code sequence have the same parity.

\section{The Map of Triangles}

\begin{figure*}
    \centering
    \begin{subfigure}{0.5\textwidth}
        \centering
        \includegraphics[scale=0.08]{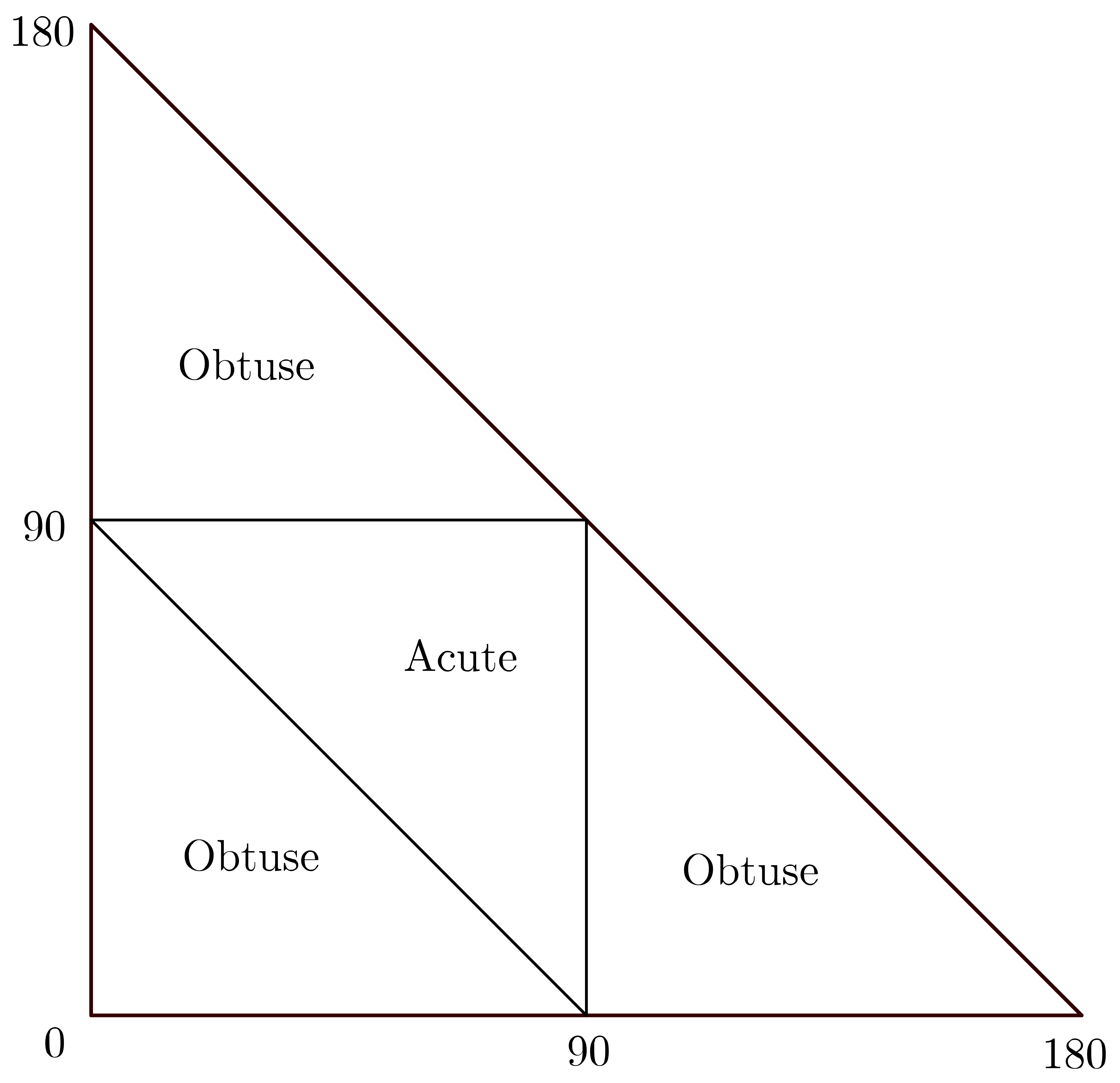}
        \caption{Acute Vs. Obtuse regions}
        \label{fig:mapa}
    \end{subfigure}%
    ~
    \begin{subfigure}{0.5\textwidth}
        \centering
        \includegraphics[scale=0.08]{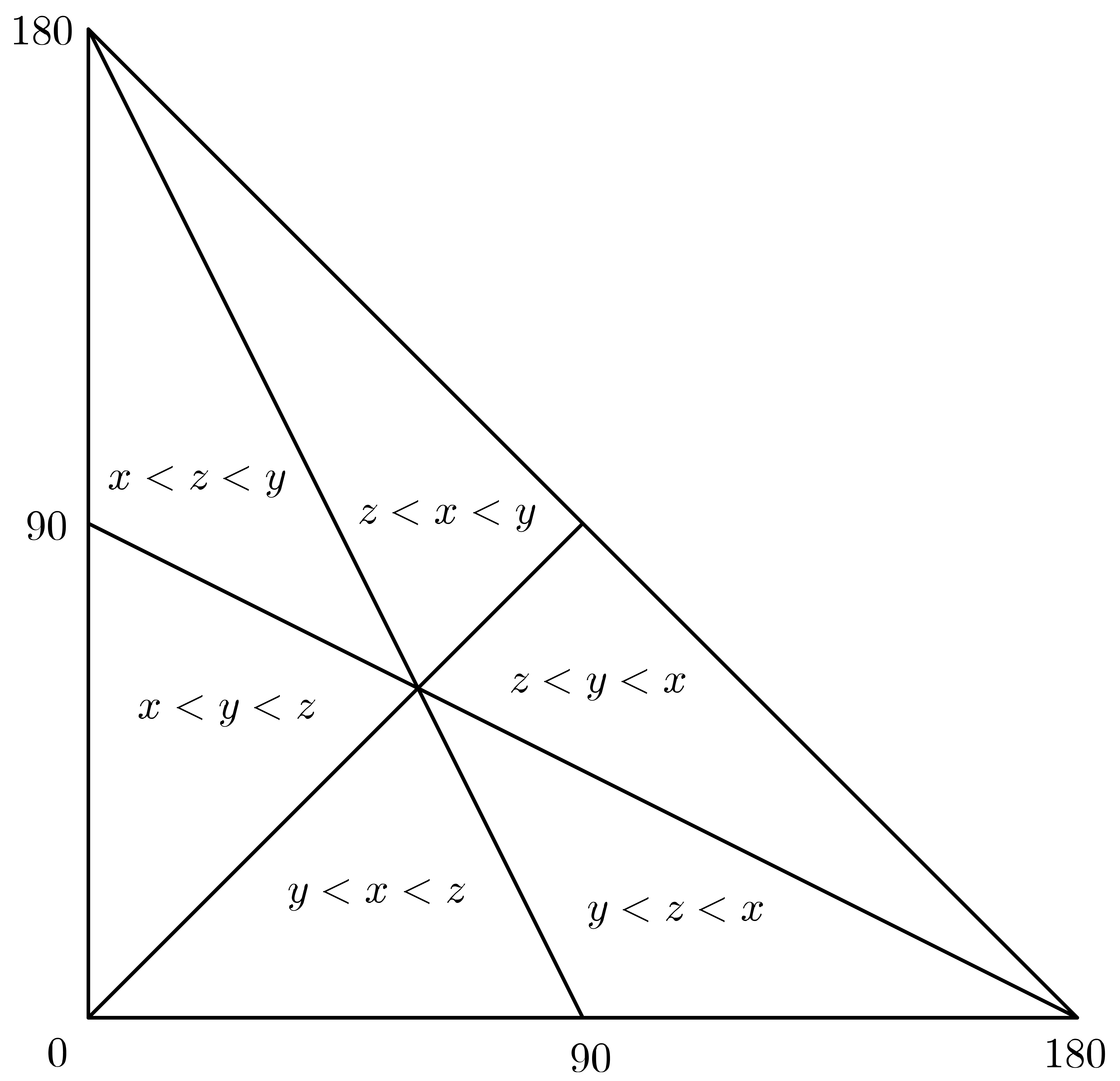}
        \caption{Angle Proportion Regions}
        \label{fig:mapb}
    \end{subfigure}
    ~
    \begin{subfigure}{0.5\textwidth}
        \centering
        \includegraphics[scale=0.08]{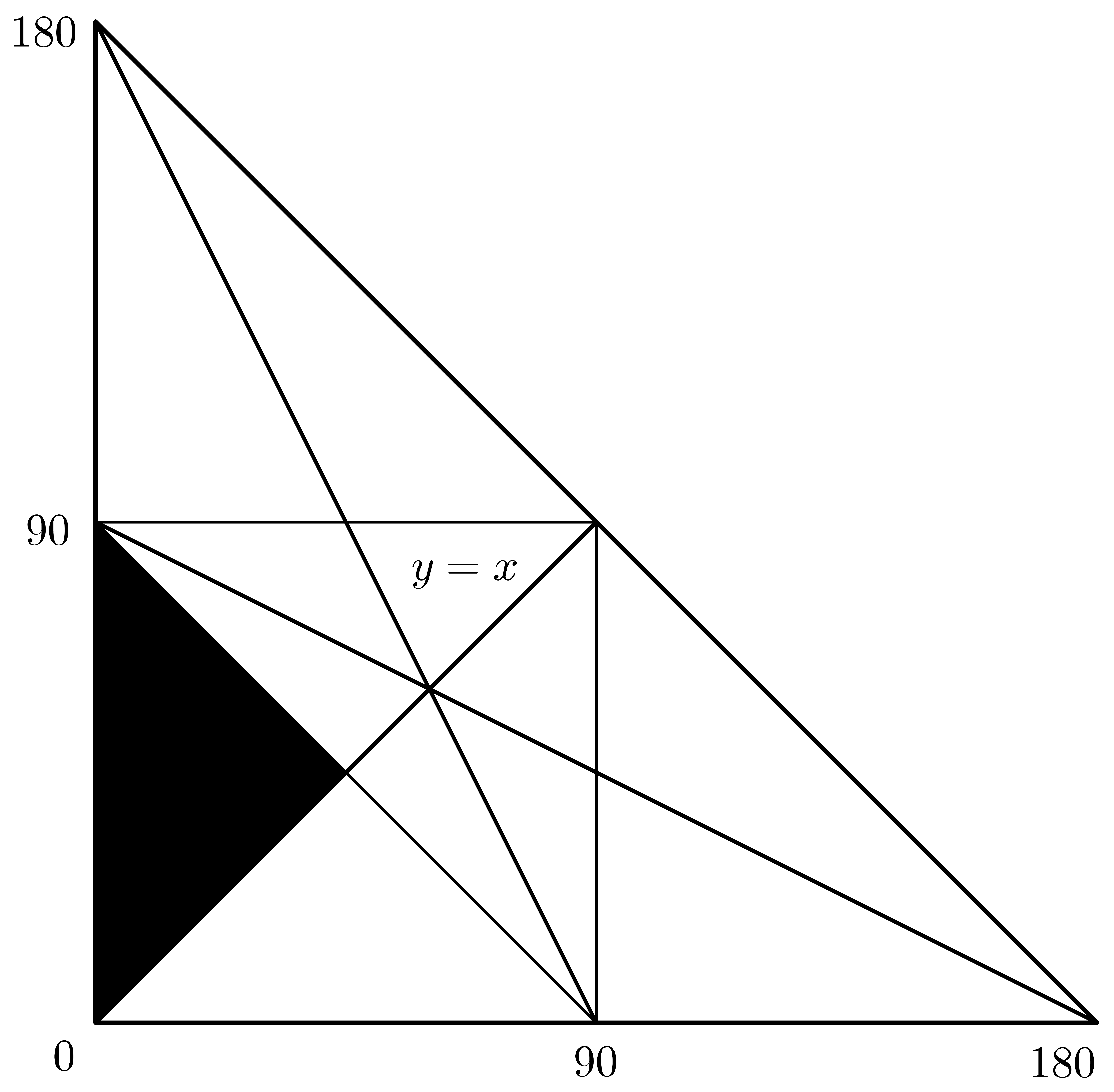}
        \caption{All Unique Obtuse Triangles}
        \label{fig:mapc}
    \end{subfigure}%
    ~
    \begin{subfigure}{0.5\textwidth}
        \centering
        \includegraphics[scale=0.16]{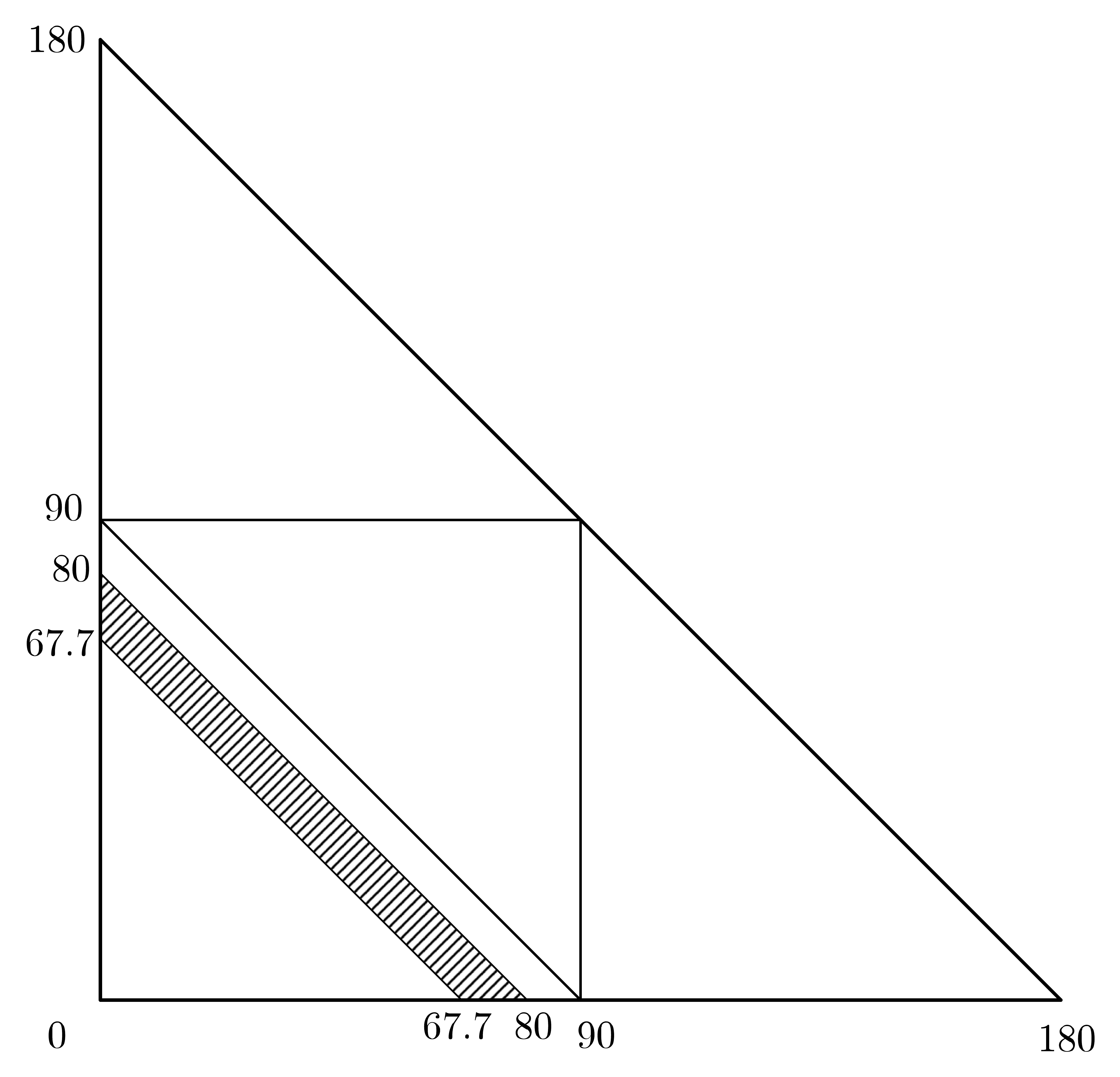}
        \caption{The 112.3 Degree Theorem}
        \label{fig:mapd}
    \end{subfigure}
    \caption{The Map of all Triangles}
\end{figure*}

We will suppose that $m<A=x, m<B=y$ and $m<C=z$ and since $z=180-x-y$, the angles of the triangle are completely determined by the coordinates $(x,y)$ and we will plot this point in an $X-Y$ coordinate system \textbf{only if that triangle has a periodic path}. Note this is independent of the size of the triangle. Observe that if the triangle corresponding to $(a,b)$ has a periodic path, then so do the triangles corresponding to $(a,180-a-b)$, $(b,a)$, $(b,180-a-b)$, $(180-a-b,a)$ and $(180-a-b,b)$. To prove that every triangle has a periodic path amounts to plotting every point $(x,y)$ with $x,y>0$ and $x+y<180$ or equivalently to plotting every point in the region $0<x\leq y\leq z$ as in Figure \ref{fig:mapb}. To prove that all obtuse triangles have a periodic path amounts to plotting every point $(x,y)$ with $x+y<90$ and $x\leq y$ which means it is enough to fill in the shaded region of Figure \ref{fig:mapc}.

Our goal in this paper is to shade in the region bounded by $x=0$, $y=0$, $x+y=80$ and $x+y=67.7$ as shown in Figure \ref{fig:mapd} or equivalently the region bounded by $x=0$, $x\leq y$, $x+y=80$ and $x+y=67.7$ which in conjunction with the other known results will prove that every triangle whose largest angle is 112.3 degrees or less has a periodic path. We call this the 112.3 degree theorem.

\section{Tower of Mirror Images of a Triangle}

\noindent

Let triangle ABC be oriented counterclockwise from A to B to C. A finite sequence of mirror images in the sides of triangle ABC will be called a \textbf{tower} if successive mirror images are in different sides of the triangle. The \textbf{length} of a tower is the number of triangles in it. We will make the convention that the first mirror image is not in side AB which will form the \textbf{base} of the tower. If the last mirror image is in side UV, then either of the other two sides can be viewed as the \textbf{top} of the tower. It is a \textbf{parallel tower} if the last mirror image in side AC or side BC makes side AB parallel to the base and pointing in the same direction as we go from from A to B. AB is then the top of the tower.

If we take a successive subset of these mirror images then we will call it a \textbf{subtower} of the given tower which is a tower in its own right allowing its base to be any one of the three sides.

If we number the triangles consecutively the first or starting triangle will always be oriented counterclockwise and so will every odd numbered triangle while every even numbered triangle will be oriented clockwise. All successive mirror images of the vertices A,B and C will be called A,B, and C points. Successive "A" points are successively labelled $A_{0}$ to $A_{n}$, successive "B" points $B_{0}$ to $B_{m}$ and successive "C" points $C_{0}$ to $C_{p}$. This means the orginal triangle is labelled $A_{0}$$B_{0}$$C_{0}$ and each vertex in the tower has an unique label.

Also note that a tower can overlap itself and that there may or may not be any poolshot associated with it as discussed in the next section. Further the A,B,C points are in fact ordered in an increasing order that follows the ordering of the formation of the sequence of mirror images of the tower. So for example, $i<j$ if and only if $C_{i}$ was formed before $C_{j}$ in the sequence of mirror images. Caution: If the tower overlaps itself, some vertices could have multiple labels which is not a problem as they would belong to different triangles.

We will color the vertices of the tower according to the following rule. $A_{0}$ is a \textbf{blue point} and $B_{0}$ is a \textbf{black point}. If the first reflection is in side AC then C is a black point and if the first reflection is in side BC then C is a blue point. Inductively if vertex U has color blue(black) and the next reflection is in side UV, then V has the opposite color black(blue).

\begin{figure}[ht]
    \centering
    \includegraphics[scale=0.15]{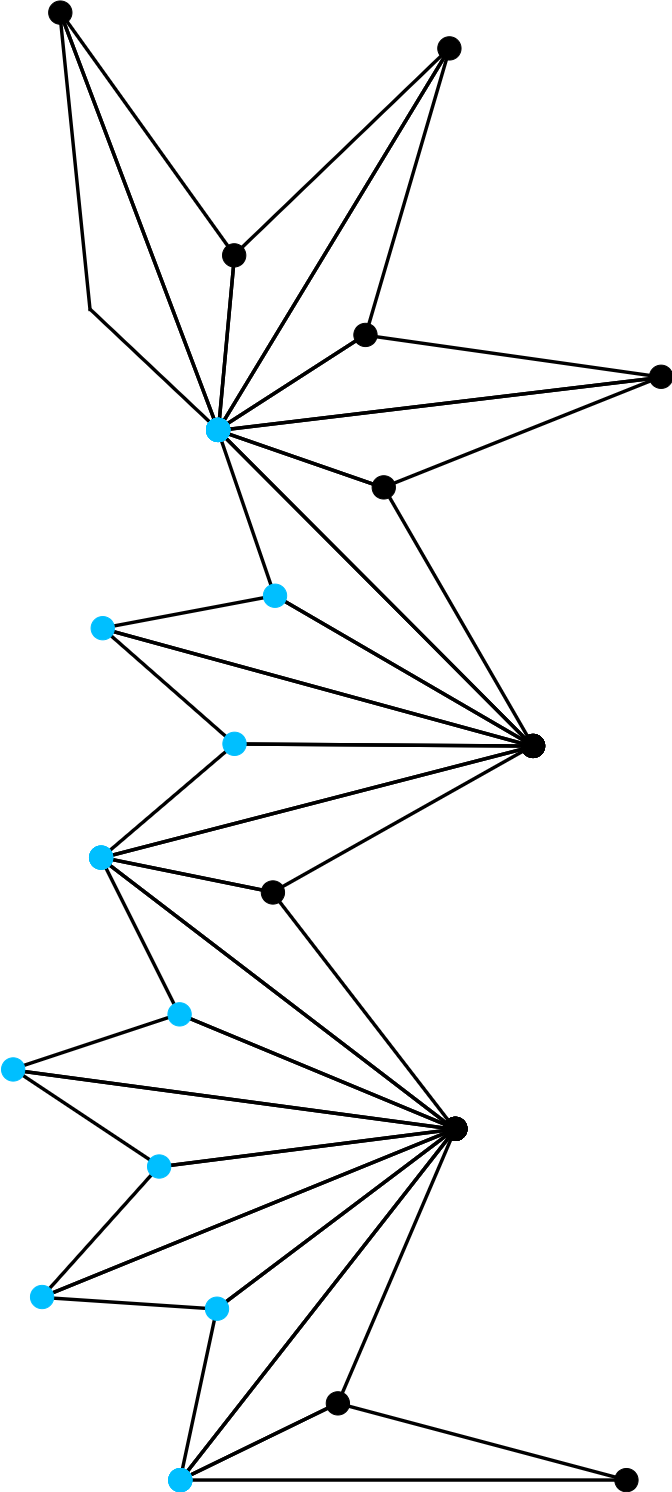}
    \caption{Tower of Mirror Images}
    \label{fig:tower}
\end{figure}

We can also use the side sequence notation to describe a tower where the first integer say 1 represents its base and successive integers represent the successive sides in which the mirror images are taken. For example the tower in figure \ref{fig:tower} can be described by the non-repeating and hence \textbf{non-legal side sequence} 131212121312121313131

\noindent \textbf{Important comment:} If a side sequence is symmetric about some integer say the integer "i" corresponding to side UV, then the corresponding subtower is symmetric about that same side UV. For example 312123\textbf{1}321213 is symmetric about the bolded 1. As a consequence, the line joining corresponding vertices of symmetric sides will be perpendicular to UV.

\section{Poolshot Towers}

Given triangle ABC oriented counterclockwise from A to B to C and given a finite billiard trajectory or poolshot starting say at side AB and which doesn't hit a vertex, if we \textbf{straighten} out this poolshot upwards by successive reflections in the sides that the poolshot hits then we get a corresponding finite \textbf{tower} of mirror images in which the poolshot is now a straight line segment. All vertices occuring on one side of the straightened poolshot will be blue points while all vertices occuring on the other side will be black points.

It follows that the convex hull of the blue points and the convex hull of the black
points are disjoint. This further means we cannot have a $blue-black-blue$ collinear situation where a black point is between two blue points. Similarly we cannot have a $black-blue-black$ collinear situation.

Since the poolshot starts at AB, if we straighten it out upwards with the base AB placed horizontal with A to the left of B and C above the base which we call \textbf{standard position}, then all blue points are on the $left$ $side$ and all black points are on the $right$ $side$ of the straightened poolshot and we get a \textbf{poolshot tower}. It is worth noting that a parallel poolshot tower must have an even number of triangles in it which is not necessarily the case for a parallel tower.

\begin{figure}[ht]
    \centering
    \includegraphics[scale=0.2]{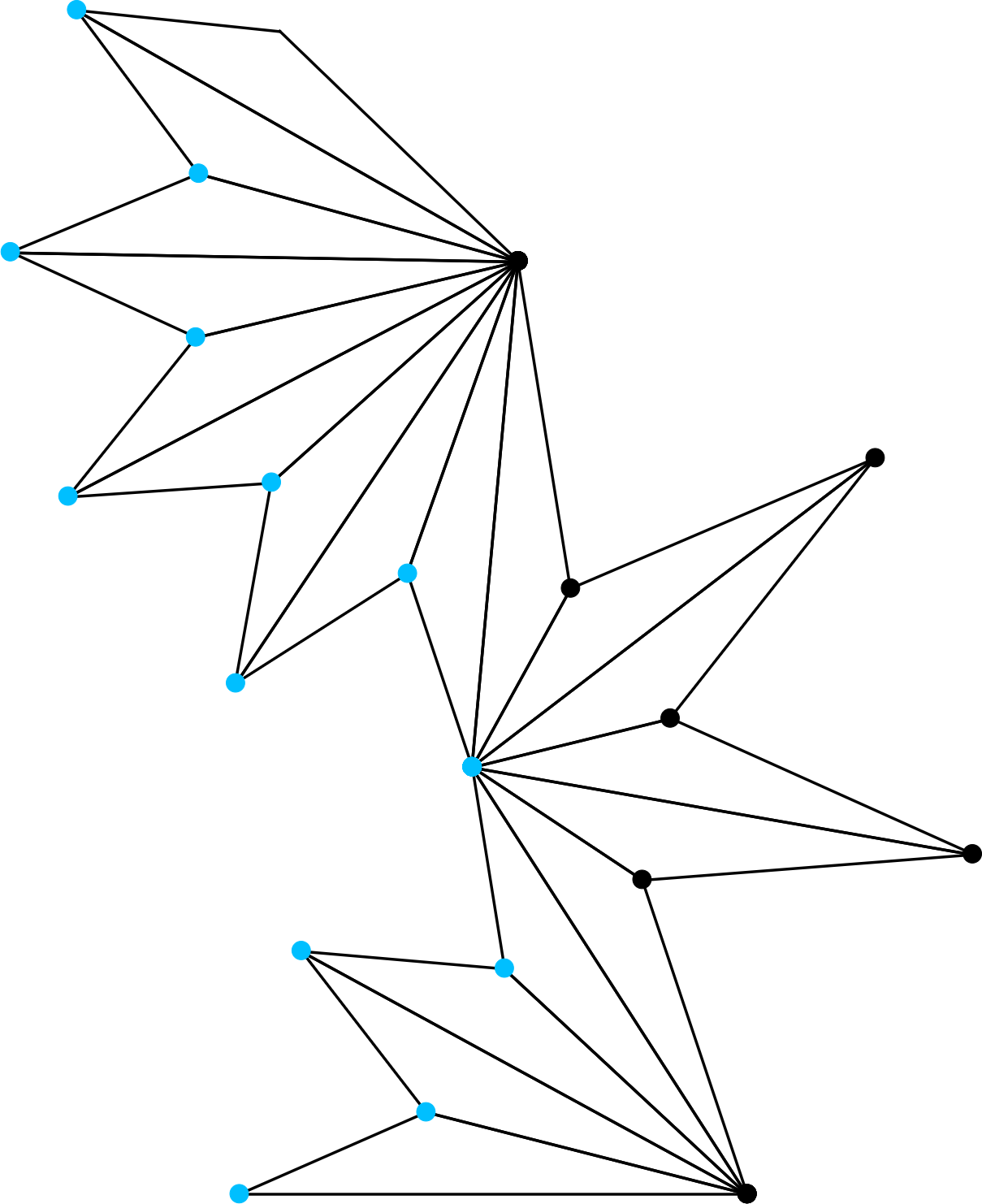}
    \caption{Poolshot Tower}
\end{figure}

\noindent
\textbf{Convention:} Any straightened poolshot can be viewed as forming the positive Y coordinate axis in an XY coordinate system by introducing a perpendicular X axis through the starting point of the poolshot on side AB. With this convention all blue points will be truly on the left side and all black points will be truly on the right side of the poolshot in this coordinate system.

\section{Periodic Poolshot Towers}

If the poolshot forms a periodic path, then the corresponding poolshot tower is called a \textbf{periodic poolshot tower}. If it starts at side $A_{0}B_{0}$ and is periodic of even length, then it finishes at $A_{n}B_{m}$ for some n and m where $A_{n}B_{m}$ is \textbf{parallel} to $A_{0}B_{0}$ where $A_{0}$ and $A_{n}$ are blue points and $B_{0}$ and $B_{m}$ are black points and we have a parallel poolshot tower. If the periodic poolshot is of odd length and finishes at $A_{n}B_{m}$, then $A_{n}B_{m}$ is \textbf{antiparallel} (in the sense that interior angles on the same side of the straightened poolshot are equal) to $A_{0}B_{0}$ with $A_{0}$ a blue point and $A_{n}$ a black point and it follows that if we double the length of the poolshot and go around the periodic path twice then $A_{2n}B_{2m}$ will be parallel to $A_{0}B_{0}$ and both $A_{0}$ and $A_{2n}$ will be on the same side of the straightened poolshot and both will be blue points and again we end with a parallel poolshot tower.

\begin{figure}[ht]
    \centering
    \includegraphics[scale=0.107]{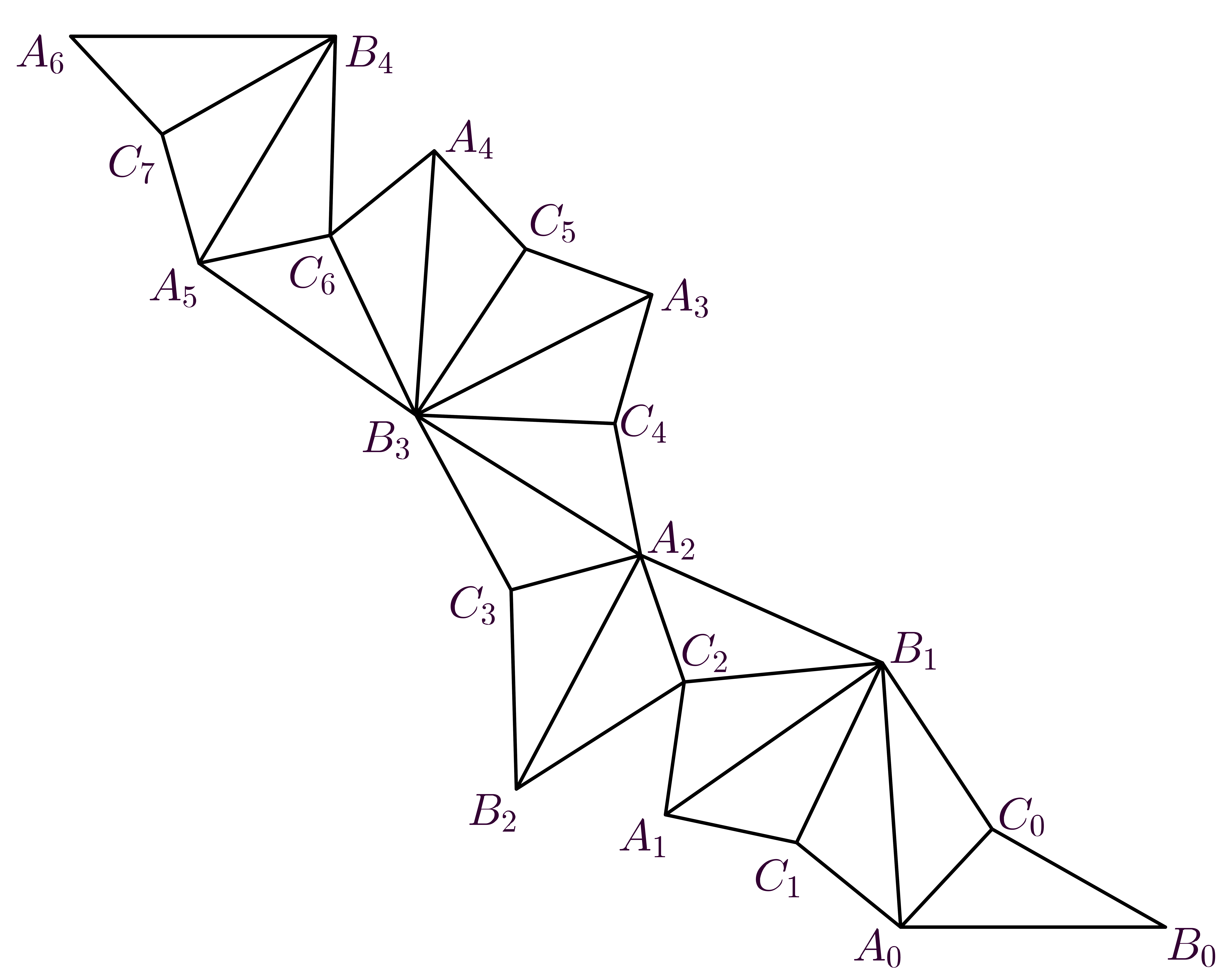}
    \caption{Periodic Poolshot Tower}
    \label{fig:periodicpoolshot}
\end{figure}

\noindent \textbf{Important comment:} In a periodic poolshot tower, the side which is at the top of the tower is completely determined and is the same as the base.

\section{The Poolshot Tower Test}

As previously noted given a poolshot tower the convex hulls of the blue and black points must be disjoint. Conversely if the convex hulls of the blue and black points respectively of a tower are disjoint then by a well known separation theorem there is a line separating the two sets and since A is a blue point and B a black point, that line must go through the base AB of the tower (and also through every segment joining a blue point to a black point) and hence there must be a straightened poolshot which produces the tower.

\textbf{The Poolshot Tower Test:} A tower is a poolshot tower if and only if the convex hulls of the blue and black points don't intersect.

\section{The Periodic Poolshot Tower Test}

If we are given a periodic poolshot tower of \textbf{even length} then as stated previously the base $A_{0}B_{0}$ and the final side $A_{n}B_{m}$ of the tower are parallel line segments. The periodic poolshot that produces the tower will leave some point $P=P_{0}$ on $A_{0}$$B_{0}$ at an angle $\theta$ where $0<\theta\leq90$ and return to that point $P=P_{q}$ on $A_{n}$$B_{m}$ also at the angle $\theta$ but on the other side of the straightened poolshot. Since $A_{0}B_{0}$ and $A_{n}B_{m}$ are parallel, this is the same acute angle between the line $A_{0}A_{n}$ and the base $A_{0}B_{0}$ or between the line $B_{0}B_{m}$ and the base $A_{0}B_{0}$. Indeed $A_{0}A_{n}$ and $B_{0}B_{m}$ are both parallel to the straightened poolshot $P_{0}P_{q}$.

Now observe that any line segment between a blue point and a black point must cross the line through $P_{0}$ and $P_{q}$. Vectorially if $v=(a,b)$ is a non-zero vector from any blue point to any black point and $w=(c,d)$ is a non-zero vector along the straightened poolshot or equivalently along $A_{0}$$A_{n}$ or $B_{0}$$B_{m}$ called the \textbf{shooting vector}, then $bc<ad$.

Conversely if we a given a parallel tower of even length in which $A_{0}B_{0}$ is parallel to $A_{n}B_{m}$ and where $bc<ad$ for every vector from a blue point to a black point then it must be a periodic poolshot tower since if U is a blue point farthest to the right (reminding the reader that we are orienting the coordinate system so that $A_{0}A_{n}$ is vertical) and V is a black point farthest to the left and since $bc<ad$, there must be a band of non zero width between the two points between which there is a periodic poolshot which produces the given tower. Hence we get the

\textbf{Periodic Poolshot Tower Test I:} A parallel tower of even length with base $A_{0}B_{0}$ parallel to $A_{n}B_{m}$ is a periodic poolshot tower if and only if $bc<ad$ (or equivalently $ad-bc>0$) for all vectors $v$ where $v=(a,b)$ is a non-zero vector from any blue point to any black point and $w=(c,d)$ is a vector from $A_{0}$ to $A_{n}$.

\noindent As with repeating side sequences we can also use a code sequence to describe a repeating tower of mirror images of triangle ABC be it a periodic poolshot tower or not. For example the periodic poolshot tower in figure \ref{fig:periodicpoolshot} can be described by the side sequence 131212313121212312 (the sequence of reflected sides) or the code sequence 2 3 1 3 5 1 1 2. Since by convention the first triangle in a tower is oriented counterclockwise, we make a similar convention as before that the first 1 in the side sequence just represents the first appearance of triangle ABC oriented counterclockwise and after that the integers represent the sequence of reflected sides and the resulting triangle that is produced. For example the last 2 in the 14th spot represents a reflection in side BC which produces the 14th triangle in the tower which must be oriented clockwise. Note that corresponding to any subcode, there is a corresponding sequence of mirror images of triangle ABC which is a subtower.

\section{Fans}

A \textbf{fan} is a tower of mirror images of a triangle in which all successive mirror images alternate between the same two sides which means that all triangles of the fan intersect at the same vertex which we will call its \textbf{center}. We can further classify the fans as \textbf{blue fans} if the center is black and all other vertices are blue points or \textbf{black fans} if its center is blue and all other vertices are black. The \textbf{central angle} of a blue or black fan is the maximal angle at its center produced by the blue or black vertices in the tower.

Any tower can be viewed as a succession of fans and in particular any poolshot tower can be viewed as a succession of alternating black and blue fans as cut off by the straightened poolshot. Even more a blue(black) fan will consist of a black(blue) center corresponding to as an example vertex A say and two circular arcs of blue(black) points on the other side of the straightened poolshot, one arc corresponding to vertex B and one arc corresponding to vertex C. The end vertices of these two \textbf{blue(black) arcs} will be called the \textbf{key blue(black) points} of this blue (black) fan. Any blue or black fan contains at most 4 key points.

\begin{figure}
    \centering

    \begin{subfigure}{0.44\textwidth}
        \centering
        \includegraphics[scale=0.11]{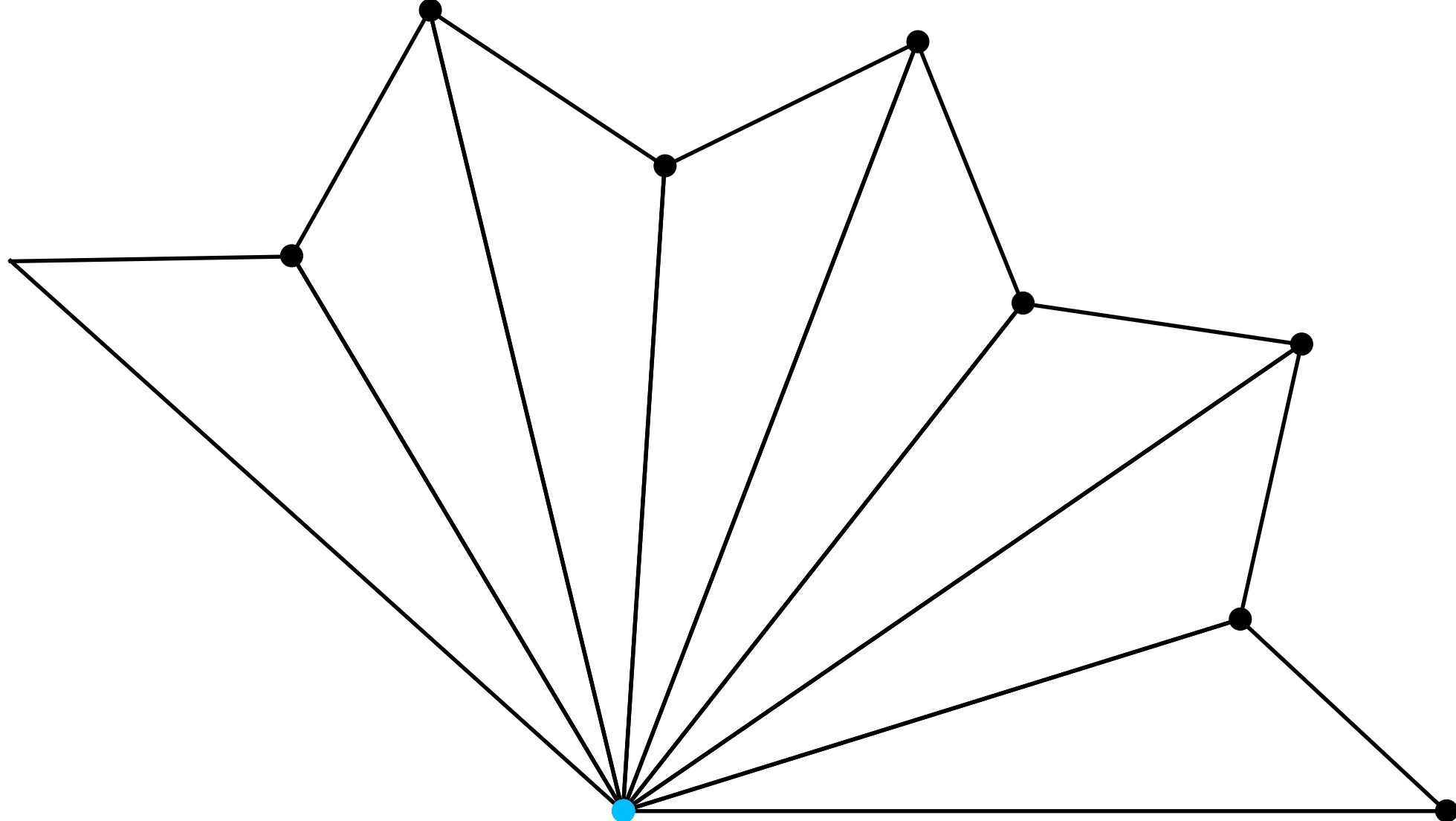}
        \caption{One black Fan}
    \end{subfigure}%
    ~
    \begin{subfigure}{0.56\textwidth}
        \centering
        \includegraphics[scale=0.13]{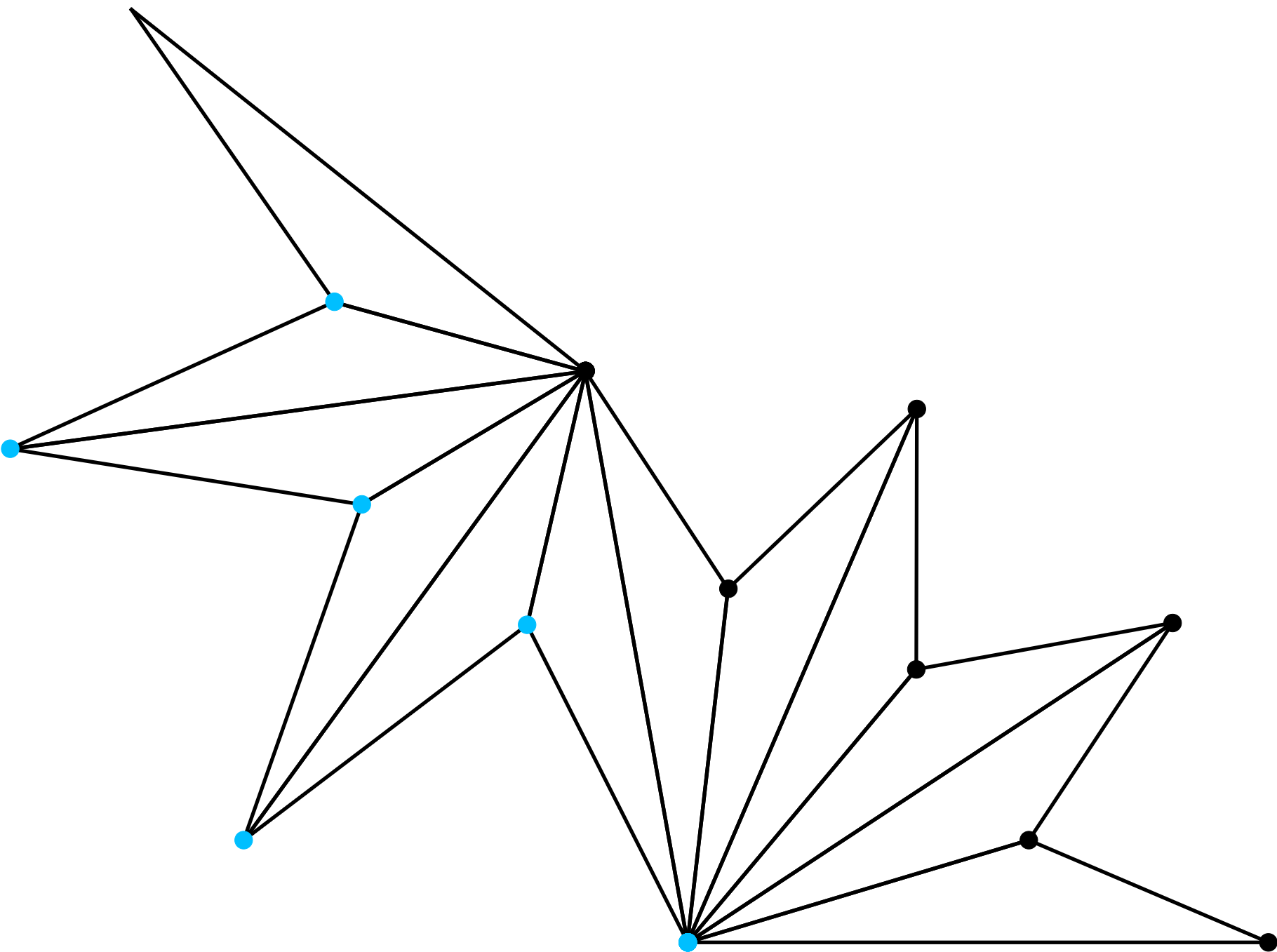}
        \caption{A Black and Blue Fan}
    \end{subfigure}

    \caption{Fan}
\end{figure}

 \textbf{Fan Fact I}: The center of a fan is a key point of the following and preceding fan of opposite color assuming it has a following or preceding fan.

\textbf{Fan Fact II}: Given a blue fan in a poolshot tower, then the key points of each blue arc are the points on the arc closest to the straightened poolshot. Similarly for black fans. This is a consequence of the fact that $sin\theta$ is a minimum at the endpoints of the interval [a,b] where $0< a \leq b <180$ or equivalently that $cos\theta$ is a minimum at the endpoints of the interval [a,b] where $-90< a \leq b <90$ and that in a poolshot tower the central angle of any fan is less than 180 degrees.

\textbf{Fan Fact III}: Every blue(black) vertex lies on some blue(black) arc whose endpoints are key blue(black) points and whose center is black(blue).

\noindent
\textbf{Labelling of the Centers and Key points in a periodic poolshot tower}

We will assume that the base is $A_{0}B_{0}$ and the first reflection is in side $A_{0}C_{0}$ , then the black point $B_{0}$ has the label $L_{(1,0)}$ and the blue point $A_{0}$ the label $L_{(2,0)}$. Now as the centers alternate between black and blue points the labels increase by one. Observe that all black centers have odd labels $L_{(2i-1,0)}$ and all blue centers have even labels $L_{(2i,0)}$ for $i\geq1$. If the tower has 2m fans in it, then the last labels are $L_{(2m+2,0)}$ and $L_{(2m+1,0)}$ which belong to the last A and B vertices in the tower respectively.

As to the key points if any belonging to the fan with center $L_{(k,0)}$, if there are four of them the first one appearing in the tower after $L_{(k-1,0)}$ is labelled $L_{(k,1)}$ and the second $L_{(k,2)}$ whereas if there is three key points, the one after $L_{(k-1,0)}$ is labelled $L_{(k,1)}$.

Note 1: The number of fans in a periodic poolshot tower in standard form equals the number of code numbers which is always even where we remind the reader that periodic paths of odd side sequence length are doubled in order to get a parallel tower.

Note 2: The first B and the last A vertices can be considered as centers of degenerate fans involving no triangles and are key points and are not counted in the number of fans.

\begin{figure}[ht]
    \centering
    \includegraphics[scale=0.7]{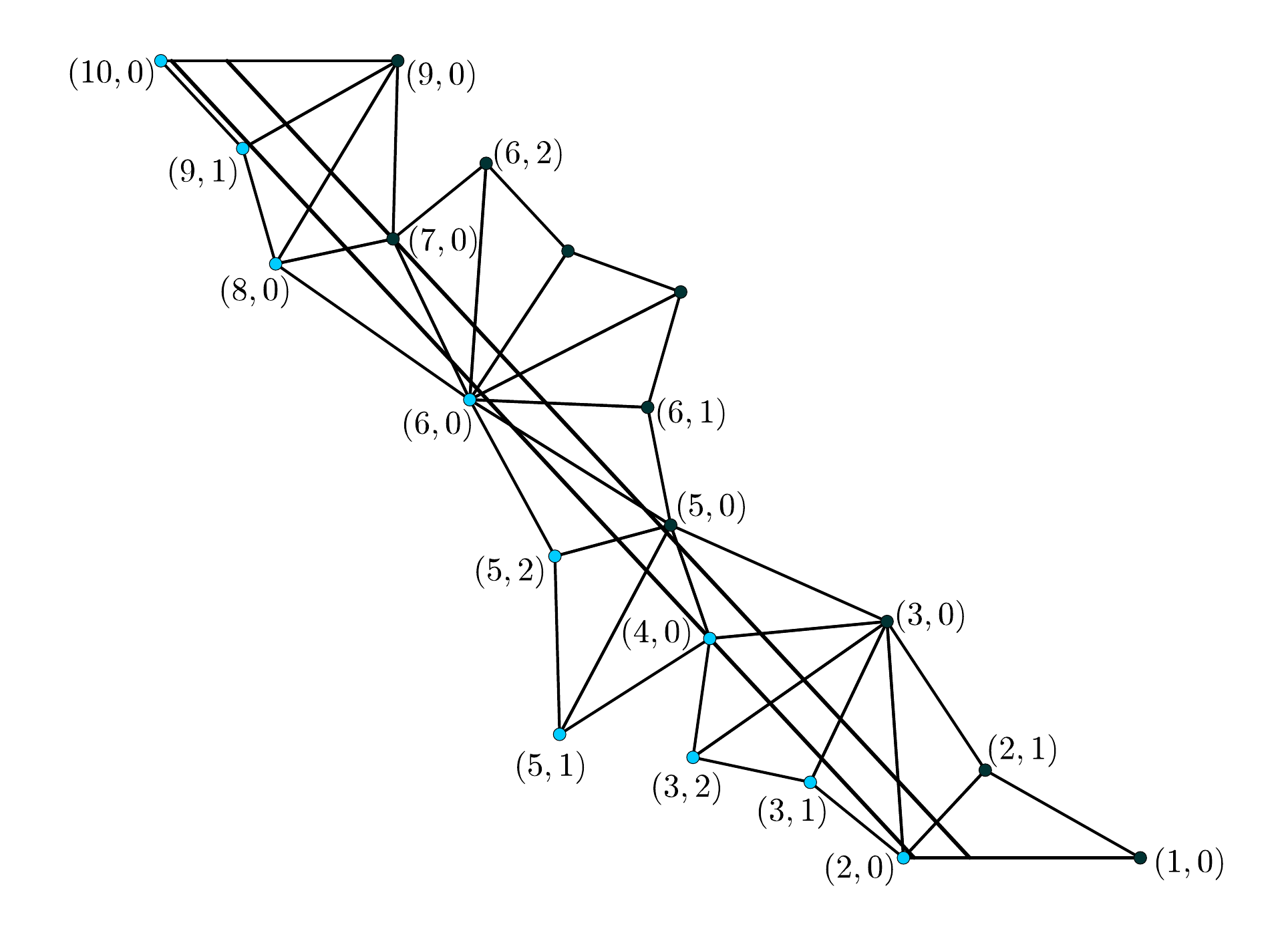}
    \caption{Labelled Periodic Poolshot Tower}
\end{figure}

\textbf{Periodic Poolshot Tower Test II:} A parallel tower of even length with base $A_{0}B_{0}$ parallel to $A_{n}B_{m}$ is a periodic poolshot tower if and only if $bc<ad$ or equivalently $ad-bc>0$ where $v=(a,b)$ is a non-zero vector from any key blue point to any key black point and $w=(c,d)$ is a vector from $A_{0}$ to $A_{n}$ \textbf{provided the central angle of any blue or black fan is less than 180 degrees}.

\begin{proof}
Choose our coordinate system so that $A_{0}A_{n}$ is vertical. Let V be a black point with the smallest X coordinate in our system and let $L_{2}$ be a vertical line through V and let U be a blue point with the largest X coordinate and let $L_{1}$ be a vertical line through U. Then our tower is a periodic poolshot tower if and only if U lies to the strict left of V.

If the tower is a periodic poolshot tower then by the Test I $ad-bc>0$ for all vectors from blue points to black points and hence this is certainly true for all vectors from key blue points to key black points.

On the other hand to show the other direction, it is enough to show that U and V are both key points since then the condition that $ad-bc>0$ will guarantee that U lies to the strict left of V. Now since V is black, it must lie on some black arc whose center is blue. If that blue center lies to the left of or on $L_{2}$, then V must be a key black point otherwise since the central angle is less than 180 degrees one of the black endpoints of the black arc through V would lie to the left of V. If that blue center lies to the right of $L_{2}$, then since that center is a key blue point and the endpoints of the black arc through V are key black points, the condition $ad-bc>0$ would force the central angle of that black arc to be greater than 180 degrees which is impossible. Hence this case can't arise and V is a key black point. Similarly U is a key blue point.
\end{proof}

Note 1: We can disregard the last blue and black points from this calculation since if $ad-bc>0$ using a vector from the first blue point $A_{0}$ to some black point (not the last), then it follows that $ad-bc>0$ using a vector from the last blue point $A_{n}$ to that same black point since both vectors are on the same side of $A_{0}A_{n}$.
Also we can disregard using a vector from $A_{n}$ to $B_{m}$, since this works if and only if the vector from $A_{0}$ to $B_{0}$ works since these are the same vectors.

Note 2: If two blue points form a vector parallel to the shooting vector then using either blue point and any fixed black point produces the same sign for $ad-bc$. This means we need only choose one blue point and disregard the others if they form a vector parallel to the shooting vector.

Note 3: If two or more blue points determine a vector parallel to the shooting vector and two or more black points determine a vector parallel to the shooting vector, then we need only choose one blue point and one black point to find the sign of $ad-bc$.

\noindent
Conclusion: In our computer calculations to show that a periodic path exists in a given triangle we use this second test taking into account the notes above.

\noindent
\textbf{Important comment:} Alternating succesive code numbers can be taken to represent alternating blue and black centers of alternating black and blue fans and we can call them \textbf{blue or black code numbers}. Observe that blue(black) code number gives one less than the number of black(blue) vertices in the corresponding fan and that the sum of the blue(black) codes plus one is the number of black(blue) vertices in the tower. Alternately each code number gives you the number of triangles in each fan and the sum of the code numbers gives you the number of triangles in the corresponding tower.

With successive code numbers in any code sequence, we can associate the X,Y or Z angles used in the central angles of the fans of the corresponding tower. Note that we use X and x, Y and y and Z and z represent the same angles. \textbf{Notationally we like to use X,Y,Z when dealing with the alternating angles of the fans in a tower and x,y,z when dealing with the angles of an individual triangle.}

Algorithm to do this:
Rule 1. Let the first code number correspond to X and the second code number correspond to Y (Note X and Y can be replaced by any of X,Y or Z)

Rule 2. Now consider any two successive code numbers $C_{i}$ and $C_{i+1}$ that have angles say X and Y assigned to them. Then if $C_{i+1}$ is even, then $C_{i+2}$ has the same angle as $C_{i}$ (X in this example) whereas if $C_{i+1}$ is odd then $C_{i+2}$ has an angle different from that of $C_{i}$ or $C_{i+1}$ (Z in this example).

Notationally we will write these angles successively above and below the corresponing code numbers starting with X say on top. We will call these the \textbf{top angles and the bottom angles}.

\noindent
Example

\noindent
X \: Z \: Y

\noindent
1 3 3 1 3 3

\noindent
\; Y \: X \; Z

Observe that the center angles of the fans are then successively 1X, 3Y, 3Z, 1X, 3Y, 3Z, and these centers alternate from one side to the other of the associated tower. These center angles are of the form code number times appropriate angle of the triangle.
Also observe that if the sum of the top angles times its code number equals the sum of the bottom angles times its code number and if the base of the tower is AB, then the side of the triangle at the top of the tower is parallel to the base although not necessarily the same side as the base.
Finally observe that if a periodic poolshot tower has an even number of triangles in it and the above is satisfied then the top and the base are parallel and of the form AB where A is a blue point and B is a black point.

Given a periodic poolshot tower with corresponding code sequence n m ... then one can determine the successive acute angles \textbf{(the shooting angles)} as the straightened poolshot crosses each side of a triangle in the tower. If the first shooting angle from the base of the tower is $\theta$ where $0<\theta<90$ and the first fan it crosses involves the angle X and has central angle nX and if n=2k+1 then the successive shooting angles are $\theta$, $\theta + x$, $\theta + 2x$, ...
, $\theta + kx$, $180-\theta -(k+1) x$, $180-\theta -(k+2) x$, ... , $180-\theta -(2k+1) x$. If n=2k, then the successive shooting angles are $\theta$, $\theta + x$, $\theta + 2x$, ...
, $\theta + (k-1)x$, $\theta + kx$//$180 - \theta - kx$, $180-\theta -(k+1) x$, $180-\theta -(k+2) x$, ... , $180-\theta -2kx$ where the // indicates that either $\theta + kx$ or $180 - \theta - kx$ is the acute angle.

Note that the first shooting angle cannot be 90 since in any fan that contains the 90 degree shooting angle, the side with the 90 degree shooting angle is also right in the center of the fan. Also observe that all the angles are integer linear combinations of x, y, 180 and $\theta$ and that as the poolshot passes from one fan to the next all the shooting angles are completely determined. In particular, if $\theta$ can be expressed as an integer linear combination of x,y and 90 then so too can every shooting angle.

\section{Classifying Codes}

\subsection{Stable and Unstable}

A code sequence is \textbf{stable} if the sum of all the top angles involving X times its code number equals the sum of all the bottom angles involving X times its code number and similarly for all the angles involving Y and all the angles involving Z.

\noindent
There are stable codes for which there are no triangles (x,y) which have a periodic path corresponding to that code sequence. 

\noindent
\textbf{Convention:} We will only call a code sequence stable if there is a triangle (x,y) which has a periodic path corresponding to that code sequence. 

If that is the case, then since the code is stable there is a finite periodic path within the triangle whose straightened trajectory is a postive minimum distance from all vertices. This means we can always change the coordinates (x,y) by a small amount and have another different triangle with the same periodic path. The upshot of this is that a region corresponding to a stable code is an open non-empty set in the plane and would cover a positive area. We will call it the \textbf{region} corresponding to that stable code sequence.

\noindent
\textbf{Convention:} We will call a code sequence \textbf{unstable} if it is not stable and the sum of the top angles times its code number equals the sum of the bottom angles times its code number and there exists a triangle (x,y) which has a periodic path corresponding to that code sequence.

For example 1 2 1 2 is an unstable code sequence since along the line 2X=2Y+2Z or X=Y+Z or X=90 or Y=90 or X+Y=90 there is a periodic path of type 1 2 1 2. In fact it can be easily shown that every right triangle has a periodic path of this type.

Observe that an unstable code sequence correponds to a \textbf{region} which is a finite open line segment whose equation is determined as above.

\subsection{The Five Code Types}

There are exactly five code types where four have even length, two of type stable or unstable and either with a 90 degree reflection or none. The last code type is of odd length which must be stable and can't have a 90 degree reflection.\newline

\noindent \textbf{CS codes:} These are stable even codes which contain a 90 degree reflection. The first such example is 1 1 1 1 2 1 1 1 1 2\newline

\noindent
Properties:

\noindent
1. Their side sequence length is a multiple of 4 and code sequence length is even and the code sequence passes the stable test.

\noindent
2. They are stable codes of the form $E_1$ $C_1$ $C_2$ ... $C_k$ $E_2$ $C_k$ ... $C_2$ $C_1$ or some cyclic permutation of the above where $E_1$ and $E_2$ are even code numbers and when converted to side sequence form becomes a side sequence of odd length where the middle number represents the side at which the pool shot hits at 90. For example if $E_1=4$ then the corresponding side sequence is of the form uvuvu and the poolshot hits the middle "u" at 90 and the rest of the shooting angles follow. Note that this means there are two special parallel sides corresponding to $E_1$ and $E_2$ in the corresponding tower which are perpendicular to the poolshot and half the length of the tower apart. We will call these the two \textbf{special perpendiculars}. It follows that these are the only sides of triangles in the tower which are perpendicular to the poolshot.

\noindent
3. The corresponding region covers an finite non zero area in the plane and is an open set in the plane.

\noindent
4. One can determine the \textbf{first} shooting angle of each fan by the following algorithm which we illustrate by the following example. Consider the CS code sequence\newline

X~ Z~~ Y\:~ Z\:\;~ X

1 1 1 1 2 1 1 1 1 2

~~Y~~X~~~X~~~Y~~~Z\newline

Then the first shooting angle of each successive fan is assuming the first one is $\theta$ are\newline

$\theta$

$180-\theta-1X$

$\theta + 1X-1Y$

$180-\theta-1X+1Y-1Z$

$\theta + 1X-1Y+1Z-1X$

$180-\theta-1X+1Y-1Z+1X-2Y$

$\theta + 1X-1Y+1Z-1X+2Y-1X$

$180-\theta-1X+1Y-1Z+1X-2Y+1X-1Z$

$\theta + 1X-1Y+1Z-1X+2Y-1X+1Z-1Y$

$180-\theta-1X+1Y-1Z+1X-2Y+1X-1Z+1Y-1X$\newline

\noindent Note that if we continued this pattern then the next shooting angle would be
$\theta + 1X-1Y+1Z-1X+2Y-1X+1Z-1Y+1X-2Z$=$\theta$ since $1X-1Y+1Z-1X+2Y-1X+1Z-1Y+1X-2Z=0$ since the code sequence is stable.

\noindent Further since there is a 90 angle associated with the first 2, the first shooting angle of that fan is $(180-2Y)/2$=90-Y and we can determine that $\theta$=X+Y-90 and hence express all shooting angles in terms of X,Y and 90.

\noindent
5. The first shooting angle of the first fan can be used to get the vector (c,d)=(-$cos\theta$, $sin\theta$) used in the periodic poolshot tower test.\newline

\begin{figure}[ht]
    \centering
    \includegraphics[scale=0.13]{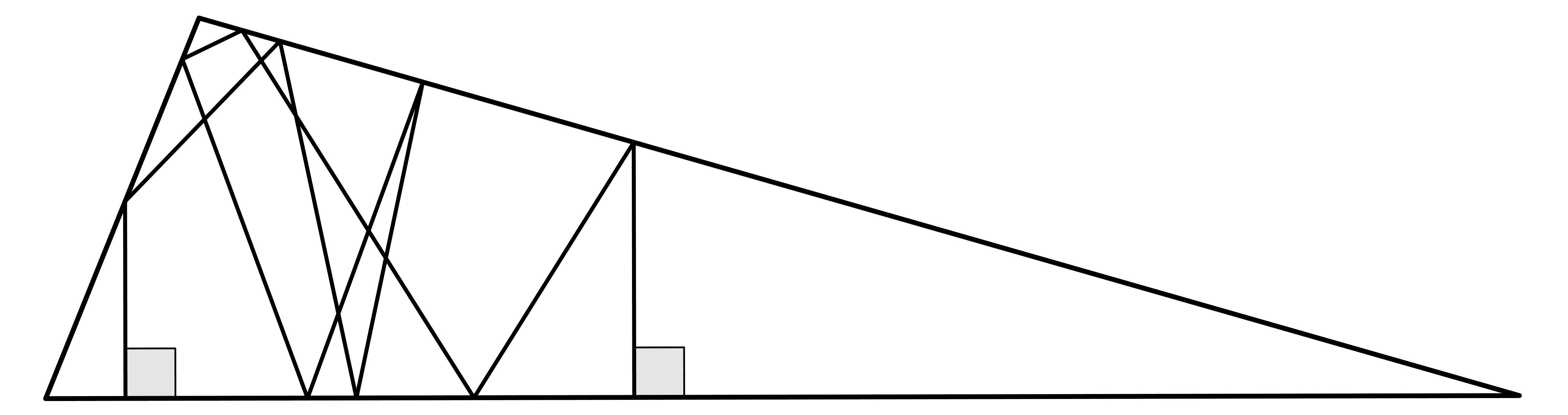}
    \caption{CS Periodic Path}
\end{figure}

\noindent \textbf{CNS codes:} These are unstable even codes which contain a 90 degree reflection. The first example is 2 2\newline

\noindent
Properties:

\noindent
1. Their side sequence length is a multiple of 2 and code sequence length is even and the code sequence doesn't pass the stable test.

\noindent
2. As above they are of the form $E_1$ $C_1$ $C_2$ ... $C_k$ $E_2$ $C_k$ ... $C_2$ $C_1$ or some cyclic permutation of the above where $E_1$ and $E_2$ are even code numbers with the difference that they are not stable codes. The corresponding tower would also contain two special perpendiculars.

\noindent
3. The corresponding region is a straight line segment with equation found by taking the sum of the top angles times their code numbers minus the sum of the bottom angles times their code numbers and setting this equal to zero. For the 2 2 CNS this becomes 2X-2Y=0 or Y=X.

\noindent
4. One can determine the first shooting angle of each fan similar to the above which we illustrate by the following example. Consider the CNS code sequence\newline

Y~~~Y

1 2 1 6

~~~Z~~~X\newline

\noindent
Then the first shooting angle of each successive fan is assuming the first one is $\theta$ are

$\theta$

$180-\theta-1Y$

$\theta + 1Y-2Z$

$180-\theta-1Y+2Z-1Y$\newline

\noindent Note that if we continued this pattern then the next shooting angle would be
$\theta + 1Y-2Z+1Y-6X$=$\theta$ since $1Y-2Z+1Y-6X$=0 or $Y=90+X$ is the equation associated with this code sequence.

\noindent Further since there is a 90 angle associated with the first 2, the first shooting angle of that fan is $(180-2Z)/2$=90-Z=X+Y-90 and we can determine that $\theta=270-X-2Y=90-3X$ and hence express all shooting angles in terms of X and 90.

\noindent
5. The first shooting angle of the first fan can be used to get the vector (c,d)=(-$cos\theta$, $sin\theta$) used in the periodic poolshot tower test.

\begin{figure}[ht]
    \centering
    \includegraphics[scale=0.12]{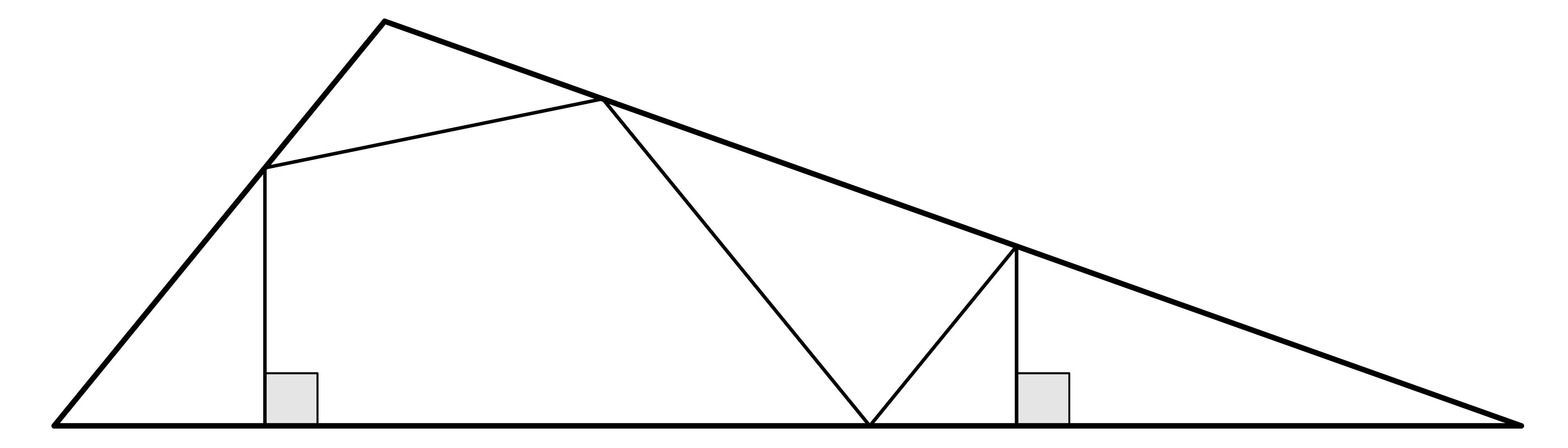}
    \caption{CNS Periodic Path}
\end{figure}

\noindent \textbf{OSO codes:} These are stable odd code sequences whose sum is an odd length or some multiple of an odd length and never contain a 90 angle. The first example is 1 1 1.\newline

\noindent
1. The side sequence length of an OSO with minimum period is odd and the code sequence length is also odd.

\noindent
2. In order to create a parallel tower from an OSO, one has to double the code as in the example above the tower has to correspond to 1 1 1 1 1 1. This means OSO's as a parallel tower are of the form $C_1$ $C_2$ ... $C_k$ $C_1$ $C_2$ ... $C_k$ or multiples thereof where $C_1$ + $C_2$ + ... + $C_k$ is odd.

\noindent
3. Since it is stable the corresponding region covers an finite non zero area in the plane and is an open set in the plane.

\noindent
4. One can determine the first shooting angle of each fan similar to the above which we illustrate by the following example. Consider the OSO code sequence. Note we don't need to double its length in this calculation.\newline

Z~~~X~~~X

1 1 2 2 5

~~~Y~~Y\newline

\noindent
Then the first shooting angle of each successive fan is assuming the first one is $\theta$ are
\newline

$\theta$

$180-\theta-1Z$

$\theta + 1Z-1Y$

$180-\theta-1Z+1Y-2X$

$\theta + 1Z-1Y+2X-2Y$\newline

\noindent Note that if we continued this pattern then the next shooting angle would be
$180-\theta-1Z+1Y-2X+2Y-5X$=$\theta$ from which we can solve for $\theta$ to get $\theta$=$(180-1Z+1Y-2X+2Y-5X)/2$.

\noindent
5. The first shooting angle of the first fan can be used to get the vector $(c,d)=(-cos\theta$, $sin\theta)$ used in the periodic poolshot tower test.

\begin{figure}[ht]
    \centering
    \includegraphics[scale=0.65]{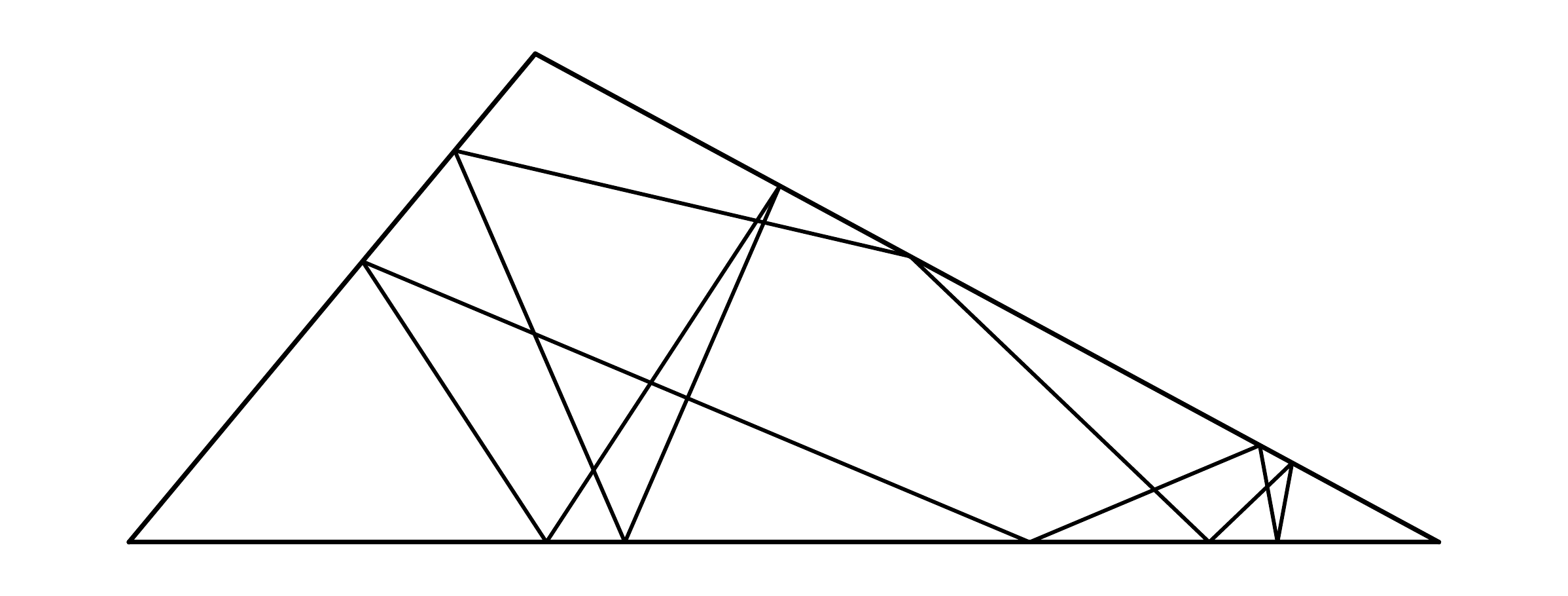}
    \caption{OSO Periodic Path}
\end{figure}

\noindent \textbf{ONS codes:} These are the unstable even codes which don't contain a 90 degree reflection. One of the first examples is 1 1 2 1 3 2.\newline

\noindent
Properties:

\noindent
1. Their side sequence length is a multiple of 2 and code sequence length is even and the code sequence doesn't pass the stable test.

\noindent
2. As above they are non stable codes of the form $C_1$ $C_2$ ... $C_k$ which are not of the CNS form.

\noindent
3. The corresponding region is a straight line segment with equation found by taking the sum of the top angles times their code numbers minus the sum of the bottom angles times their code numbers and setting this equal to zero. For the 1 1 2 1 3 2 ONS this becomes 4X-2Y=0 or Y=2X.

\noindent
4. One can determine the first shooting angle of each fan similar to the above which we illustrate by the following example. Consider the ONS code sequence\newline

X~~~X

 2 2 4 4

~~~Y~~Y\newline

\noindent
Then the first shooting angle of each successive fan is assuming the first one is $\theta$ are\newline

$\theta$

$180-\theta-2X$

$\theta + 2X-2Y$

$180-\theta-2X+2Y-4X$\newline

\noindent Note that if we continued this pattern then the next shooting angle would be
$\theta+2X-2Y+4X-4Y$=$\theta$ and we \textbf{cannot} solve for $\theta$ since on the corresponding linear tile Y=X and then $2X-2Y+4X-4Y$=0. It can be shown that the shooting angles are not integer or even rational linear combinations of X,Y and 90.

\noindent
5. Because of the above, the vector (c,d)=(-$cos\theta$, $sin\theta$) used in the periodic poolshot tower test is calculated another way as in the next section.

\noindent 6. There is a way to eliminate $\theta$ to form a bounding polygon which is discussed later and which contains the straight line segment region corresponding to the ONS code.

\begin{figure}[ht]
    \centering
    \includegraphics[scale=0.68]{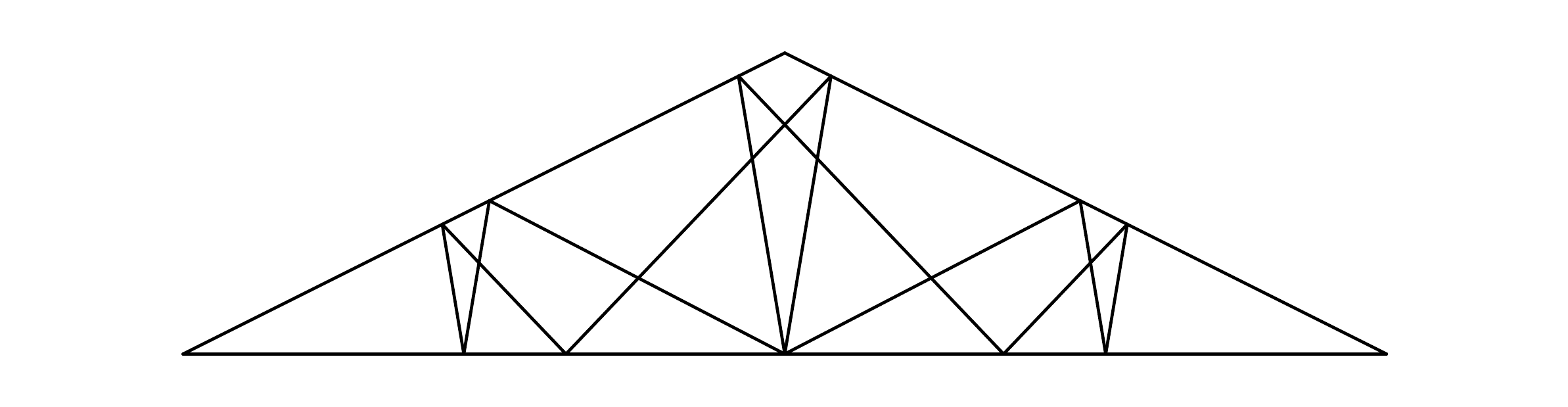}
    \caption{ONS Periodic Path}
\end{figure}

\noindent \textbf{OSNO codes:} These are the stable even codes which don't contain a 90 degree reflection. The first example is 1 1 2 2 1 1 3 3.\newline

\noindent
Properties:

\noindent
1. Their side sequence length is a multiple of 2 and code sequence length is even and the code sequence passes the stable test.

\noindent
2. They are stable codes of the form $C_1$ $C_2$ ... $C_k$ which are not of the CS or even multiples of the OSO form.

\noindent
3. The corresponding region covers an finite non zero area in the plane and is an open set in the plane.

\noindent
4. One can determine the first shooting angle of each fan similar to the above which we illustrate by the following example. Consider the OSNO code sequence\newline

X~~~Y~~~Y~~Z

 1 1 2 2 1 1 3 3

~~~Z~~~Z~~~X~~~Y\newline

\noindent
Then the first shooting angle of each successive fan is assuming the first one is $\theta$ are\newline

$\theta$

$180-\theta-1X$

$\theta + 1X-1Z$

$180-\theta-1X+1Z-2Y$

$\theta + 1X-1Z+2Y-2Z$

$180-\theta-1X+1Z-2Y+2Z-1Y$

$\theta + 1X-1Z+2Y-2Z+1Y-1X$

$180-\theta-1X+1Z-2Y+2Z-1Y+1X-3Z$\newline

\noindent Note that if we continued this pattern then the next shooting angle would be
$\theta + 1X-1Z+2Y-2Z+1Y-1X+3Z-3Y$=$\theta$ since $1X-1Z+2Y-2Z+1Y-1X+3Z-3Y=0$ since the code sequence is stable.

As above it can be shown that the shooting angles are not integer or even rational linear combinations of X,Y and 90.

\noindent
5. Because of the above, the vector $(c,d)=(-cos\theta$, $sin\theta$) used in the periodic poolshot tower test is calculated another way as in the next section.

\noindent 6. There is a way to eliminate $\theta$ to form a bounding polygon which is discussed later and which contains the open region corresponding to the OSNO code.

\begin{figure}[ht]
    \centering
    \includegraphics[scale=0.6]{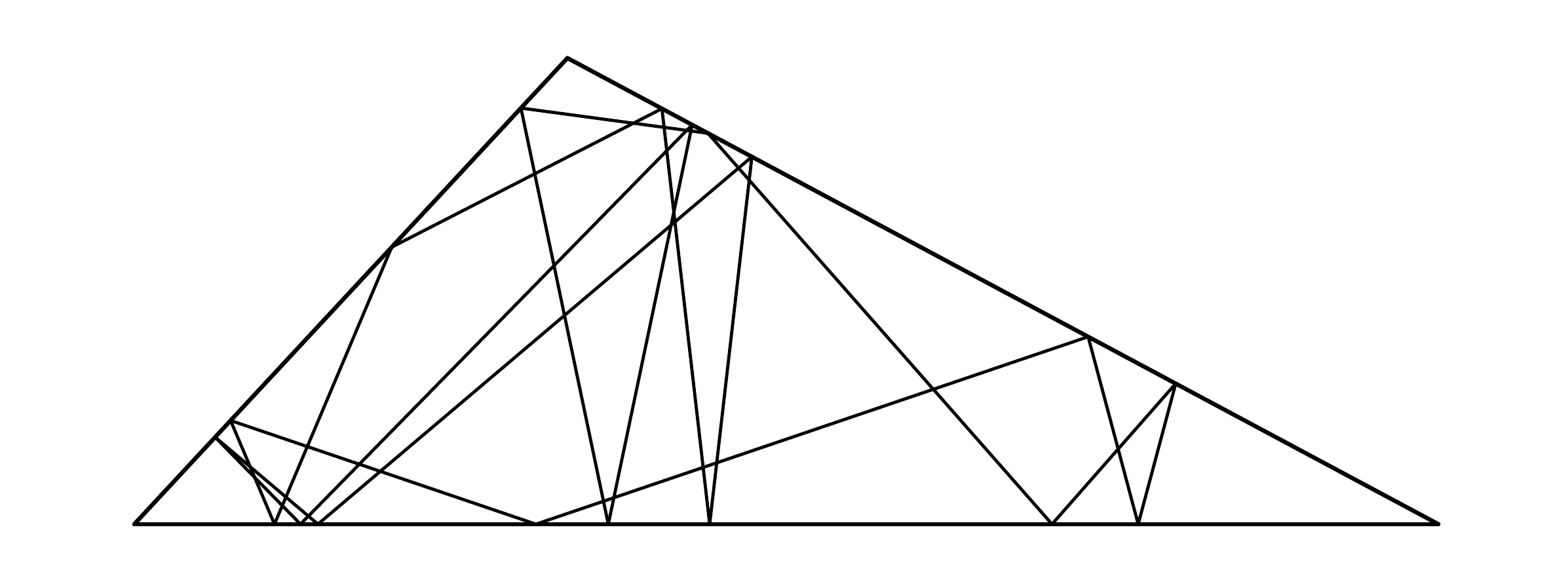}
    \caption{OSNO Periodic Path}
\end{figure}

\noindent There are faster ways to test a parallel tower for being a periodic poolshot tower if it is of the form $E_1$ $C_1$ $C_2$ ... $C_k$ $E_2$ $C_k$ ... $C_2$ $C_1$ or some cyclic permutation of the above where $E_1$ and $E_2$ are even code numbers and the sum of the top angles times its code number equals the sum of the bottom angles times its code number. This means the two special sides corresponding to $E_1$ and $E_2$ are parallel and perpendicular to $A_{0}A_{n}$.

\textbf{Periodic Poolshot Tower Test III:} A parallel tower of code sequence form $E_1$ $C_1$ $C_2$ ... $C_k$ $E_2$ $C_k$ ... $C_2$ $C_1$ (of CS or CNS form) or some cyclic permutation of the above where $E_1$ and $E_2$ are even code numbers and the sum of the top angles times its code number equals the sum of the bottom angles times its code number and with base $A_{0}B_{0}$ parallel to $A_{n}B_{m}$ and where $UV$ and $WZ$ are the two special perpendiculars associated with $A_{0}A_{n}$ is a periodic poolshot tower if and only if $bc<ad$ or equivalently $ad-bc>0$ where $v=(a,b)$ is a non-zero vector from any key blue point to any key black point which lie between or on the two special perpendiculars and $w=(c,d)$ is a vector from $A_{0}$ to $A_{n}$ \textbf{provided the angle of any blue or black fan is less than 180 degrees}.

\noindent
\section{Previous results on periodic paths:}

\noindent
1. Acute triangles always have a periodic path of OSO code type 1 1 1 which has length 3 namely the orthic triangle whose vertices are the feet of the altitudes.

\begin{figure}[ht]
    \centering
    \includegraphics[scale=0.27]{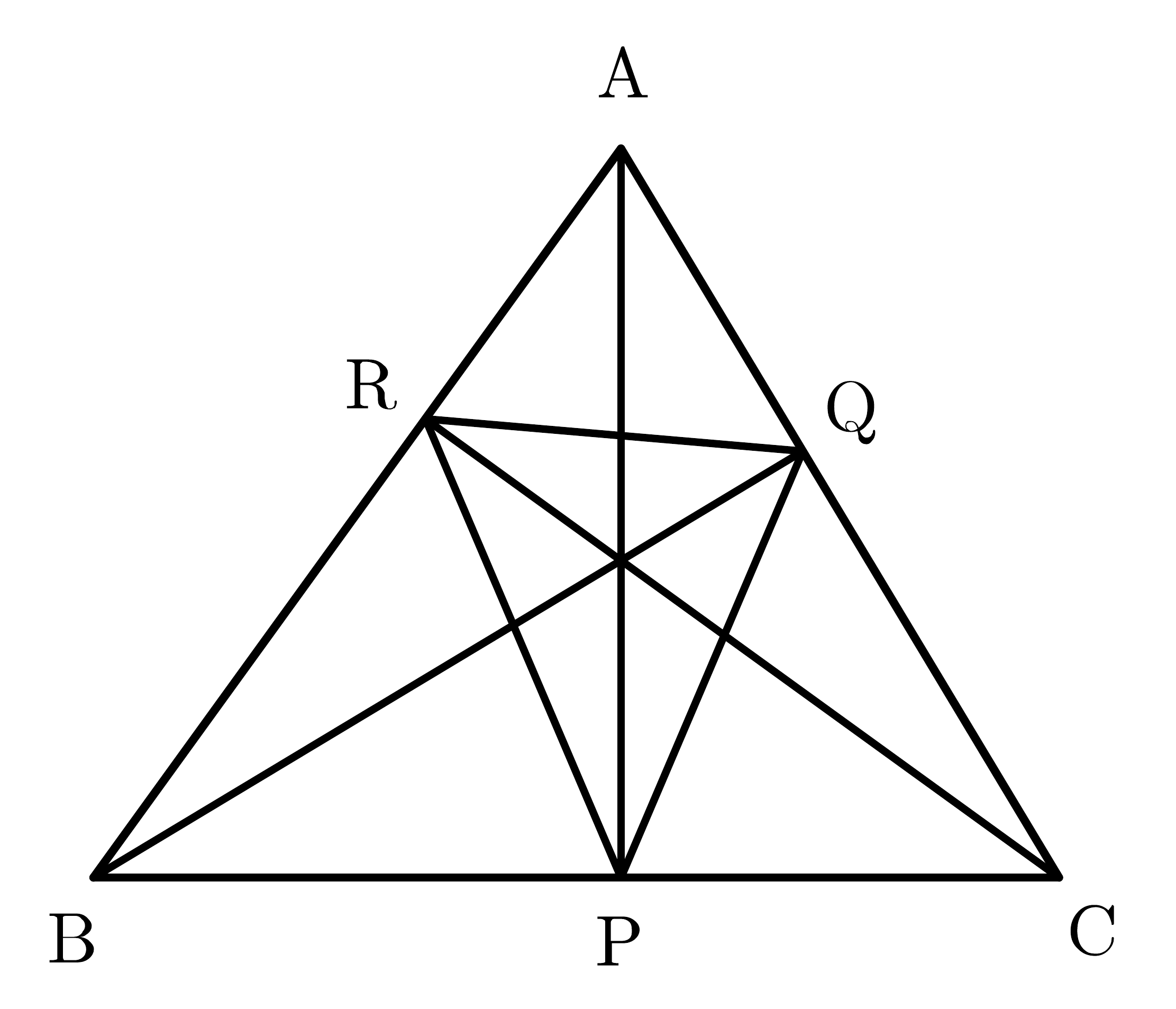}
    \caption{Orthic Triangle}
\end{figure}

\noindent
2. Right triangles always have a periodic path of CNS code type 1 2 1 2.

\begin{figure}[ht]
    \centering
    \includegraphics[scale=0.12]{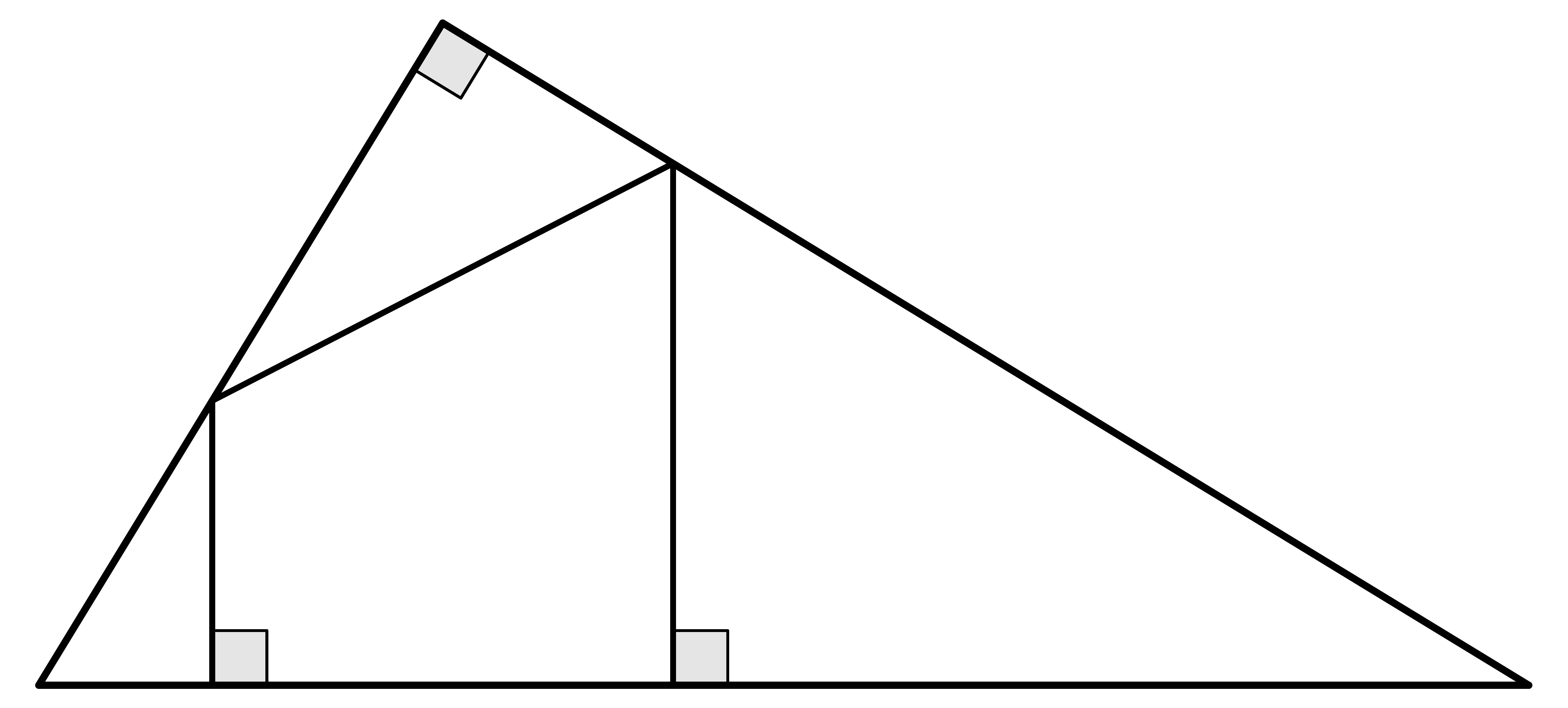}
    \caption{Right Triangle with CNS Periodic Path 1 2 1 2}
\end{figure}

\noindent
3. Isosceles triangles always have a periodic path of CNS code type 2 2. 

\begin{figure}[ht]
    \centering
    \includegraphics[scale=0.12]{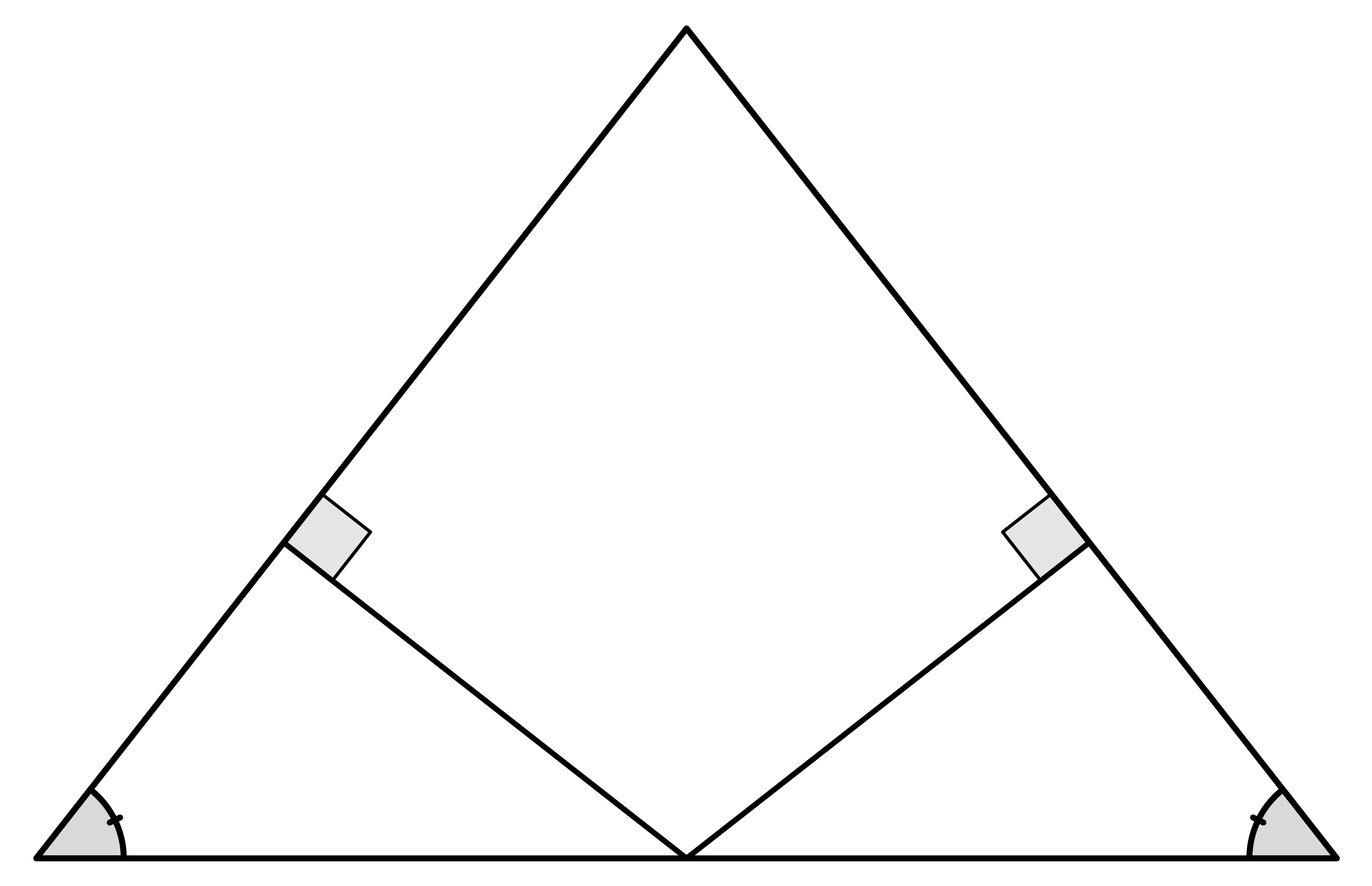}
    \caption{Isosceles Triangle with CNS Periodic Path 2 2}
\end{figure}

\noindent
4. Rational triangles always have a periodic path of type CNS. A triangle is \textbf{rational} if all angles can be expressed as integer multiples of x where x divides 90, say $m<A=nx$, $m<B=mx$ and $m<C=px$. We can prove this from the following facts about rational triangles.

\noindent
Fact I: If a poolshot leaves a side at 90 degrees and $90=qx$ for some integer $q$, then it bounces off any side at some integer multiple of x.
Conclusion: There are at most $q$ angles involved as a poolshot bounces off the sides of the triangle.

\noindent
Fact II: If a poolshot leaves a side at 90 degrees and $90=qx$ for some integer x and hits a vertex say vertex B, then it must enter B along a ray making an angle tx where $0<t<m$.

Conclusion: There are at most n+m+p-3 perpendicular poolshots which hit a vertex in a rational triangle. This means if a 90 degree poolshot leaves a side and doesn't hit a vertex, there is an open band around that poolshot of finite width $\delta>0$ which leaves that side at 90 and never hits a vertex and in which $\delta$ is maximal. This is only possible if the band hits another side at 90 and becomes periodic since otherwise the band must hit some side say side AB infinitely often at some angle jx for a fixed integer j less than q and which is repeated infintely often. But each time this band hits this side, it hits on some open interval $(a_i,b_i)$ of width $\delta$ no two of which can intersect without contradicticting the maximality of $\delta$. But since AB is finite this is impossible.

\noindent

\section{Calculating the coordinates of the vertices in a code tower}

Given a code tower of even length, let us assume that the ordering of the angles in the code tower is such that the first top angle is $U_{1}=X$ and that the last bottom angle is $U_{2k}=Z$ as shown.\newline

$U_{1}$~~~~~~~$U_{3}$

$C_1$ $C_2$ $C_3$... $C_{2k}$

~~~~$U_{2}$~~~~~~~~~~$U_{2k}$\newline

Then the base of the tower is AB with $m<A=x$, $m<B=z$ and $m<C=y$ and if we let $AB=siny$, $BC=sinx$ and $AC=sinz=sin(x+y)$ since $z=180-x-y$ we then can recursively calculate the coordinates of each fan center $L_{(i,0)}$ as follows: Letting $a_{i}=U_{i}C_i$ be the center angle of the fan with corresponding label $L_{(i+1,0)}$ for $i\geq 1$, $u_{i}$ be the length of the side between $(x_{i},y_{i})$ and $(x_{i+1},y_{i+1})$ where $(x_{i},y_{i})$ are the coordinates of the center of the fan at $L_{(i,0)}$ then recursively let $x_{1}=sin y$, $y_{1}=0$, $x_{2}=0$, $y_{2}=0$ and
\newline

\noindent$x_{2n}=x_{2n-1}-u_{2n-1}cos(a_{2n-2}-a_{2n-3}+a_{2n-4}...-a_{3}+a_{2}-a_{1})$

\noindent$y_{2n}=y_{2n-1}+u_{2n-1}sin(a_{2n-2}-a_{2n-3}+a_{2n-4}...-a_{3}+a_{2}-a_{1})$

\noindent$x_{2n+1}=x_{2n}+u_{2n}cos(a_{2n-1}-a_{2n-2}+a_{2n-3}...+a_{3}-a_{2}+a_{1})$

\noindent$y_{2n+1}=y_{2n}+u_{2n}sin(a_{2n-1}-a_{2n-2}+a_{2n-3}...+a_{3}-a_{2}+a_{1})$
\newline
\newline
\noindent
An example of this process is worked through in the appendix A.

\section{The Prover}

The basic idea behind the prover is the Mean Value Theorem in two dimensions.\newline

\noindent
\textbf{The Mean Value Theorem:} Let $f(x,y)$ be a differentiable function of two variables, then
$f(b_{1}, b_{2}) - f(a_{1}, a_{2}) = f_x(c_{1}, c_{2})(b_{1} - a_{1}) + f_y(c_{1}, c_{2})(b_{2} - a_{2})$ for some
($c_{1}, c_{2}$) on the line between ($a_{1}, a_{2}$) and ($b_{1}, b_{2}$).\newline

\noindent
\textbf{Fact I:} If $f(x,y)=\sum_{i=1}^{k} \pm u_icos(m_{i}x + n_{i}y)$, then $f_x=\sum_{i=1}^{k} \mp$$m_{i}u_{i}$sin($m_{i}x+ n_{i}y$) and $f_y$=$\sum_{i=1}^{k} \mp$$n_{i}u_{i}$sin($m_{i}x$ + $n_{i}y$) and hence
$\lvert f_x\rvert\leq\sum_{i=1}^{k} $$\lvert m_{i}u_{i}\rvert=M$ and $\lvert f_y\rvert\leq\sum_{i=1}^{k} $$\lvert n_{i}u_{i}\rvert=N$

\noindent
Similarly if $f(x,y)=\sum_{i=1}^{k} \pm u_isin(m_{i}x + n_{i}y)$\newline

\noindent
\textbf{Fact II:} If $f(b_{1}, b_{2}) >0$ and $f(a_{1}, a_{2}) \leq0$, then by the Mean Value Theorem
$f(b_{1}, b_{2}) - f(a_{1}, a_{2}) = f_x(c_{1}, c_{2})(b_{1} - a_{1}) + f_y(c_{1}, c_{2})(b_{2} - a_{2})$$\geq$$f(b_{1}, b_{2})$
and since $ f_x(c_{1}, c_{2})(b_{1} - a_{1}) + f_y(c_{1}, c_{2})(b_{2} - a_{2})\leq M\lvert b_{1} - a_{1}\rvert + N\lvert b_{2} - a_{2}\rvert$, then $f(b_{1}, b_{2})$$\leq$$M\lvert b_{1} - a_{1}\rvert + N\lvert b_{2} - a_{2}\rvert$\newline

\noindent
\textbf{Conclusion:} If $f(b_{1}, b_{2})> M\lvert b_{1} - a_{1}\rvert + N\lvert b_{2} - a_{2}\rvert$, then $f(a_{1}, a_{2}) >0$\newline

\noindent
\textbf{Fact III:} Let $(b_{1}, b_{2})$ be the center of a square of side $2r>0$, then for any $(a_{1}, a_{2})$ in or on the boundary of the square we must have $(b_{1} - a_{1}) \leq r$ and $(b_{2} - a_{2}) \leq r$ and hence
$M\lvert b_{1} - a_{1}\rvert + N\lvert b_{2} - a_{2}\rvert \leq (M+N)r$\newline

\noindent
\textbf{Conclusion I:} If $f(b_{1}, b_{2}) > (M+N)r \geq M\lvert b_{1} - a_{1}\rvert + N\lvert b_{2} - a_{2}\rvert$
then $f(a_{1}, a_{2})>0$\newline

\noindent
\textbf{Conclusion II:} If $0<r<f(b_{1}, b_{2})/(M+N)$, and $f(b_{1}, b_{2})>0$, then all points $(a_{1}, a_{2})$ in the square centered at $(b_{1}, b_{2})$ and of side 2r must also satisfy $f(a_{1}, a_{2})>0$, \newline

\noindent
\textbf{The Gradient Algorithm:} For every function $f_j(x,y)=\sum_{i=1}^{k} \pm u_icos(m_{i}x + n_{i}y)$ or

\noindent$f_j(x,y)=\sum_{i=1}^{k} \pm u_isin(m_{i}x + n_{i}y)$ and which form the boundary equations of a stable code region\newline

\noindent
1. calculate $G_j=\sum_{i=1}^{k} \lvert u_i\rvert (\lvert m_i\rvert + \lvert n_i\rvert )$

\noindent
2. then if a square centered at $(b_1,b_2)$ satisfies $f_j(b_1,b_2)>0$ and
$f_j(b_1,b_2) - rG_j >0$ for all $f_j$ where 2r is the length of a side of the square, then every point in or on the boundary of the square satisfies $f_j>0$ and that square lies completely within the given code region.\newline

Note: If we are dealing with a CNS or ONS code and its corresponding linear region, then it is exactly the same algorithm as long as $(b_1,b_2)$ is the center of the linear region intersecting the square.

\noindent
\textbf{The Triple Rule:} Suppose two code regions $R_1$ and $R_2$ intersect along a common boundary line segment from $(x_1,y_1)$ to $(x_2,y_2)$. This can happen if $R_1$ is defined by the equations $f_i(x,y)>0$ and $f(x,y)>0$ and $R_2$ is defined by the equations $g_i(x,y)>0$ and $g(x,y)>0$ and the equations $f(x,y)=0$ and $g(x,y)=0$ have a common factor of the form $sin(ax+by)$ or $cos(ax+by)$. Now observe that $sin(ax+by)=0$ if and only if $ax+by=180k$ and $cos(ax+by)=0$ if and only if $ax+by=90+180k$ for some integer k.

\begin{figure}[ht]
    \centering
    \includegraphics[scale=0.5]{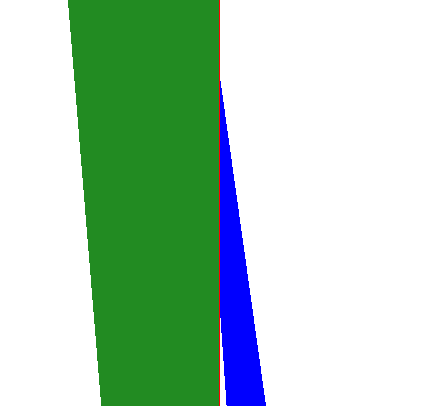}
    \caption{Two Stable Regions Sharing a Line Segment}
\end{figure}

\noindent Lets illustrate with the sine case and let us further suppose that $R_1$ lies between the parallel lines $ax+by=180(k-1)$ and $ax+by=180k$ and that $R_2$ lies between the parallel lines $ax+by=180k$ and $ax+by=180(k+1)$ and that $sin(ax+by)>0$ between the first two parallels and $sin(ax+by)<0$ between the second two parallels. Let $f(x,y)=sin(ax+by)u(x,y)$ and $g(x,y)=sin(ax+by)v(x,y)$

\noindent Now consider a square with sides parallel to the coordinate axis whose vertex coordinates are all rational numbers and which lies between the first and third parallels and which may or may not intersect the second parallel. It is worth noting that if this square lies inside either code region then it cannot intersect any of the three parallel lines above since $sin(ax+by)$ is zero there. We can use the following to decide if every point in the square including its boundary has a periodic path.

\noindent
\textbf{Triple Rule Algorithm:} Using interval arithmetic, if we can prove \newline

1. that each of the four corners $(x,y)$ of the square satisfy $180(k+1)>ax+by>180(k-1)$

2. that the center of the square $(x_0,y_0)$ satisfies $f_i(x,y)>0$, $g_i(x,y)>0$, $u(x,y)>0$ and $-v(x,y>0$ (noting that we use -v since we assume $sin(ax+by)<0$ between the second pair of parallels) and each of these equations satisfies the Gradient algorithm with respect to the given square. It then follows that all points on the square and its boundary satisfy these inequalities.

3. that each point on the common boundary line segment $ax+by=180k$ from $(x_1,y_1)$ to $(x_2,y_2)$ has a periodic path corresponding to a CNS or ONS code which includes this line segment and runs from from $(x_3,y_3)$ to $(x_4,y_4)$ and that each of the four corners of the given square lies between the lines $x=x_3$ and $x=x_4$ or between the lines $y=y_3$ and $y=y_4$.\newline

Then that square must lie within $R_1$ union $R_2$ union $R_3$ where $R_3$ is the linear region corresponding to the CNS or ONS code from 3 and every point in that union has a periodic path.

\noindent Proof: If all points on the square satisfy $sin(ax+by)>0$, it is within $R_1$. If all points satisfy $sin(ax+by)<0$, it is within $R_2$. Otherwise it is within $R_1$ union $R_2$ union $R_3$.

\noindent
QED

\section{The Two Infinite Patterns}

\noindent There are two infinite patterns which converge to the line segment x=0, $67.5<y<90$ as follows.

\begin{figure}[ht]
	\centering
    \begin{subfigure}{0.5\textwidth}
		\centering
		\includegraphics[scale=0.12]{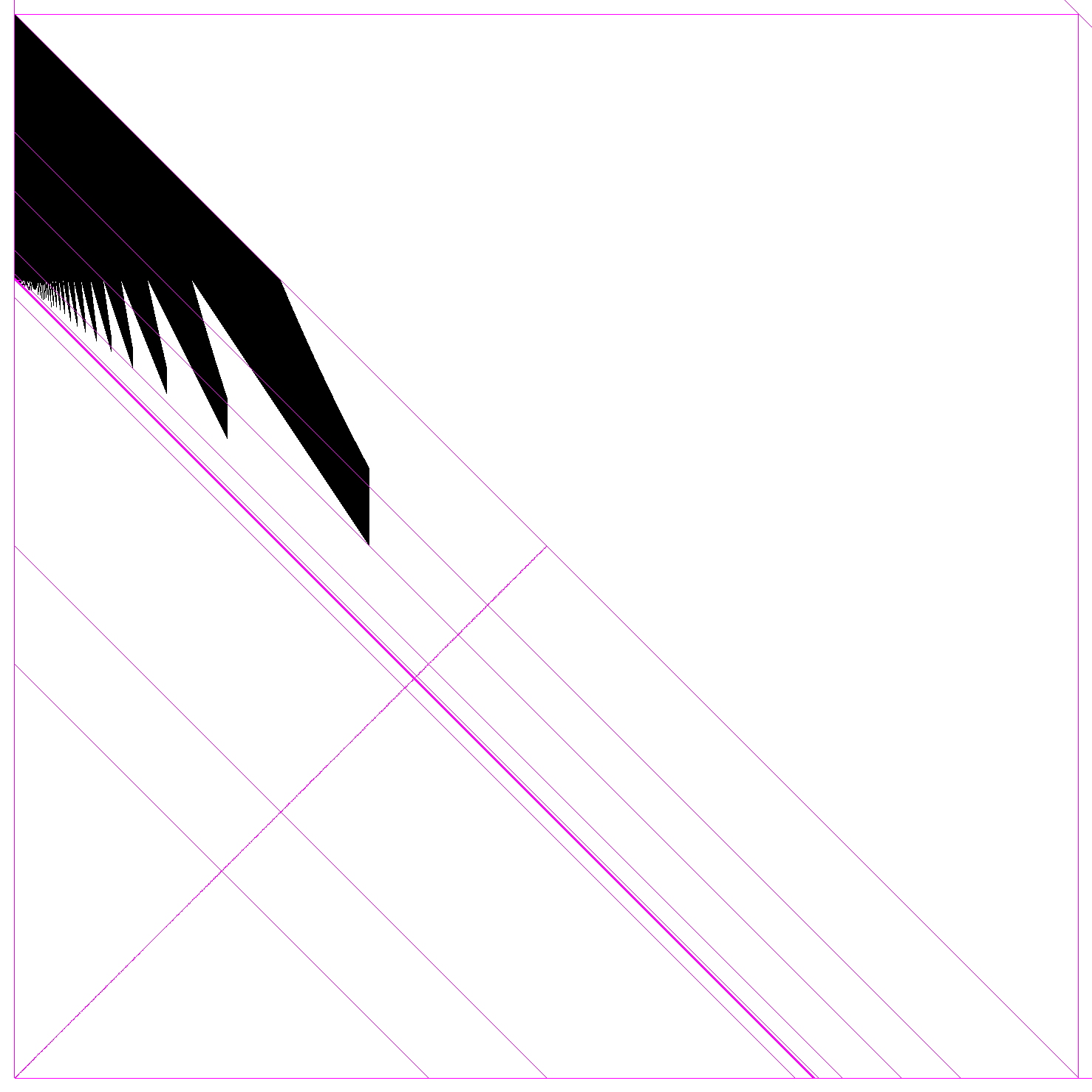}
		\caption{Infinite Pattern 1}
	\end{subfigure}%
	~
    \begin{subfigure}{0.5\textwidth}
		\centering
		\includegraphics[scale=0.12]{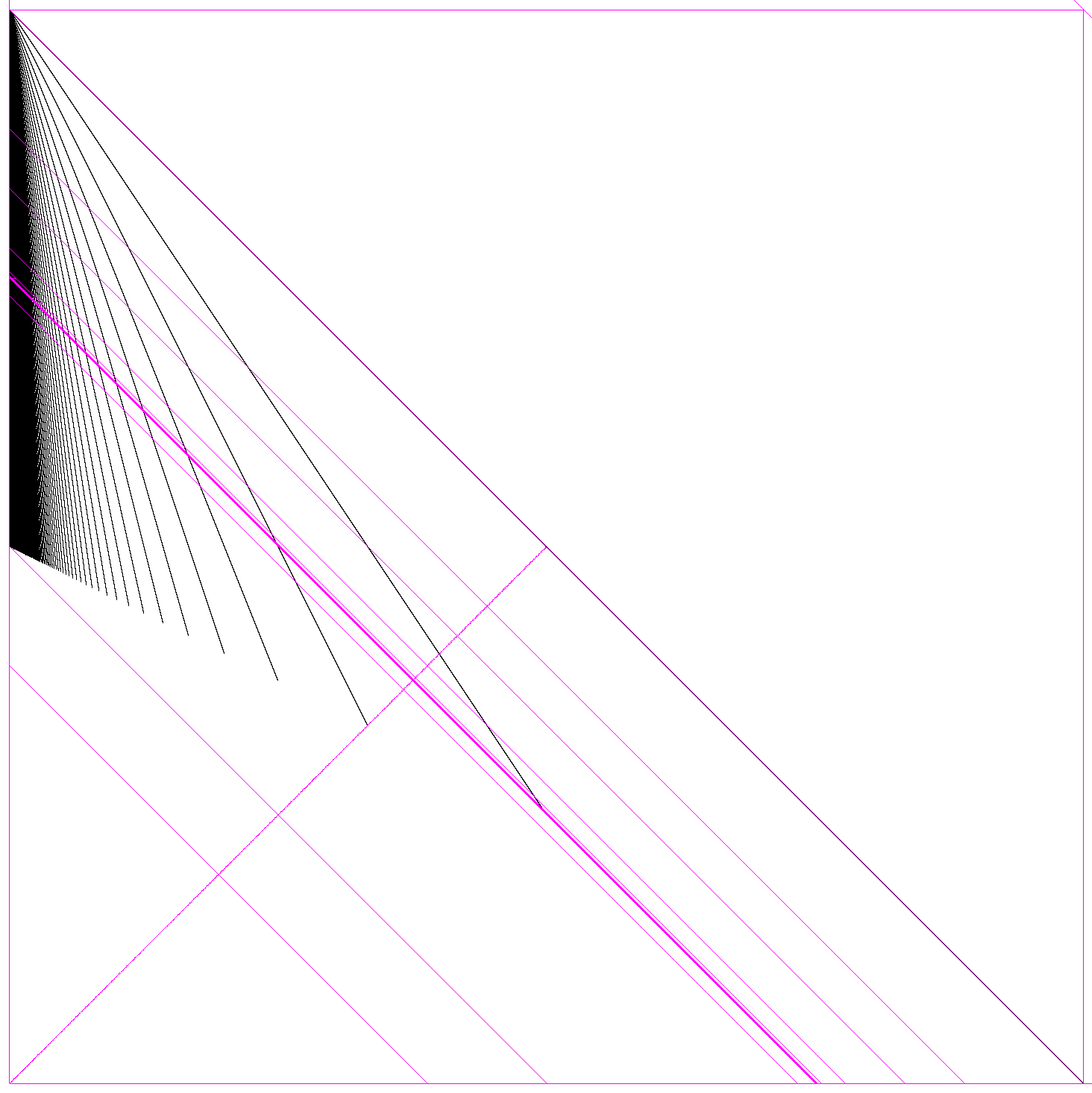}
		\caption{Infinite Pattern 2}
	\end{subfigure}
	\caption{The two infinite patterns}
\end{figure}\newpage

\textbf{Infinite Pattern I:}
Given a triangle ABC with $m<A=x$, $m<B=y$ where $(n+1)x+2y<180<(n+2)x+2y$ and $0<x<90/(2n+2)$ for $n\geq1$, then it contains a CS periodic path 1 1 $2n+1$ 1 2 1 $2n+1$ 1 1 $4n+2$. Note: Since $x<22.5$ and $y+(n+2)x/2>90$ and $(n+2)x/2<(n+2)22.5/(2n+2)<22.5$ then $y>67.5$. Also observe that since $n\geq1$, then $2x+2y\leq(n+1)x+2y<180$ which means that $x+y<90$. Finally observe that the successive regions determined by these conditions share the boundary line $(n+1)x+2y=180$ for $n\geq2$.

\begin{figure}[ht]
    \centering

    \begin{subfigure}{1\textwidth}
        \centering
        \includegraphics[scale=0.18]{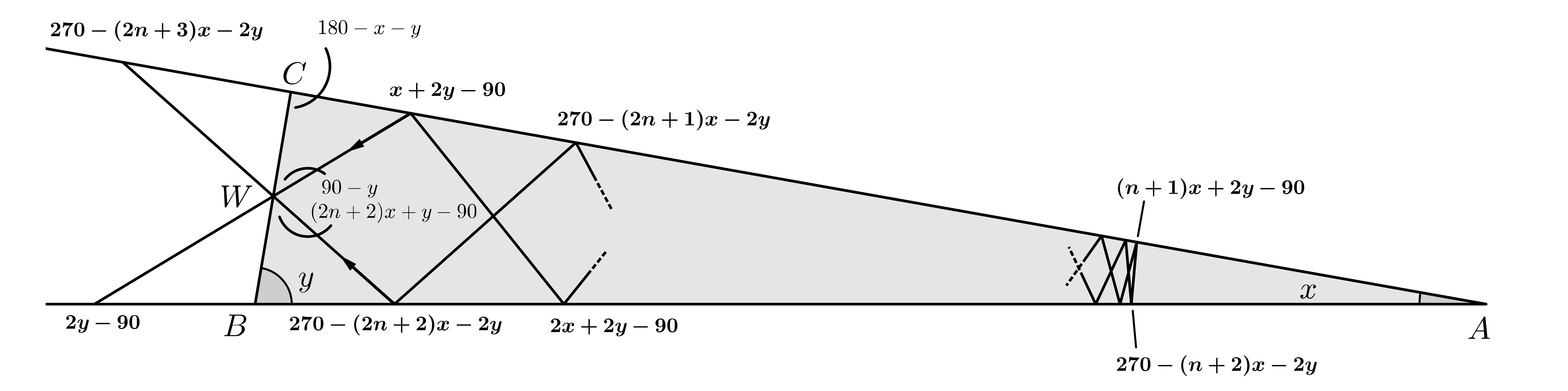}
        \caption{$n\geq2$ even case}
        \label{fig:infa}
    \end{subfigure}
    \begin{subfigure}{1\textwidth}
        \centering
        \includegraphics[scale=0.18]{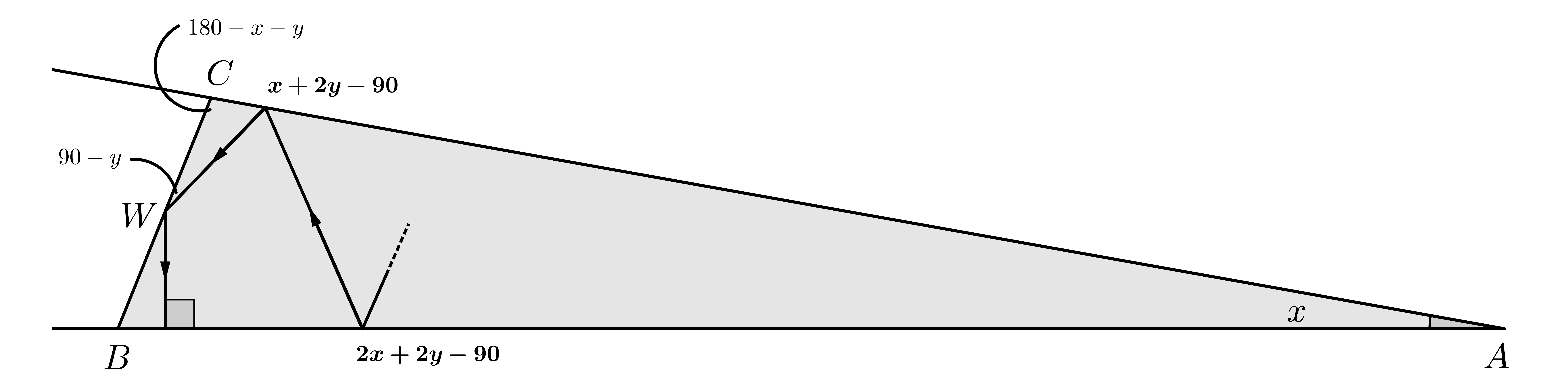}
        \caption{$n\geq2$ even case}
    \end{subfigure}
    \begin{subfigure}{1\textwidth}
        \centering
        \includegraphics[scale=0.18]{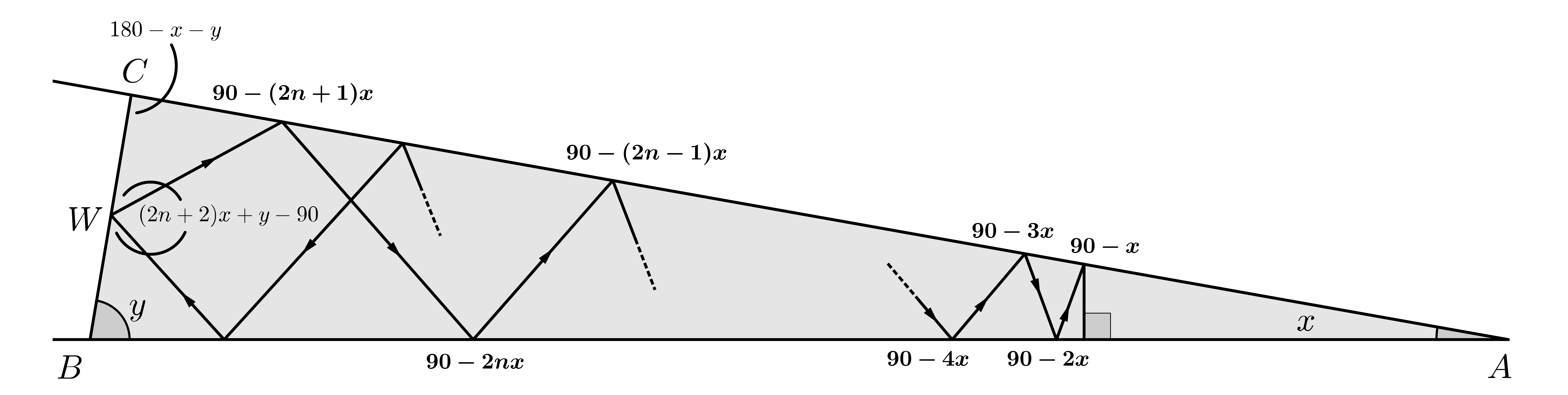}
        \caption{$n\geq2$ even case}
    \end{subfigure}
	
    \caption{Infinite Pattern I}
\end{figure}
\begin{proof}
For the $n\geq2$ even case, take a wedge of acute angle x and vertex A with arms $l_1$ and $l_2$ as shown on figure \ref{fig:infa} and pick a point P on $l_1$ and shoot a poolball at an acute angle $0<270-(n+2)x -2y<90$ as shown. Now since $270-(n+2)x -2y+x<180$, it hits $l_2$ at an angle $0<(n+1)x+2y-90<90$ and continues bouncing all on the wedge at angles $nx+2y-90>(n-1)x+2y-90>...>x+2y-90>2y-90>0$ as shown. If we consider the ray from P in the other direction, it bounces off the sides at the angles $270-(n+3)x-2y>270-(n+4)x -2y>...>270-(2n+3)x -2y>0$ where the last inequality holds since $(n+1)x+2y<180$ and $(n+2)x<(2n+2)x<90$. Now let W be the intersection of the last two rays and draw a line through W hitting $l_1$ at B, $l_2$ at C and such that the angle at B as shown is y. Observe that B lies between the last two reflections on $l_1$ since $y>2y-90$ since $y>45$ and $180-y>270-(2n+2)x -2y$ since $(2n+2)x+y>90$ where this last inequality holds since $(2n+2)x+y>(n+2)x/2+y>90$. The last ray from $l_2$ hits BC at W at the angle $90-y$ and bounces to hit AB at 90 whereas the last ray from $l_1$ hits W at $(2n+2)x+y-90$ and since $(2n+2)x+y-90$ +$180-x-y$ = $90+(2n+1)x<90+(2n+2)x<180$, it reflects off BC and hits $l_2$ at $90-(2n+1)x$. Observing that $90-(2n+1)x>x$ since $(2n+2)x<90$, the ray then bounces off $l_1$ and $l_2$ until it hits at 90 producing a CS periodic path 1 1 $2n+1$ 1 2 1 $2n+1$ 1 1 $4n+2$.

The odd case is handled similarly interchanging $l_1$ and $l_2$.

\end{proof}

\textbf{Infinite Pattern II:}
Given a triangle ABC with $m<A=x$, $m<B=y$ where $(n+1)x+2y=180$ and $0<x<90/n$ for $n\geq1$, then it contains a CNS periodic path 1 2 1 2$n$.

\noindent
Note: These are just the boundary lines (extended) between the regions of theorem 1.

\begin{figure}[ht]
    \centering

    \begin{subfigure}{1\textwidth}
        \centering
        \includegraphics[scale=0.18]{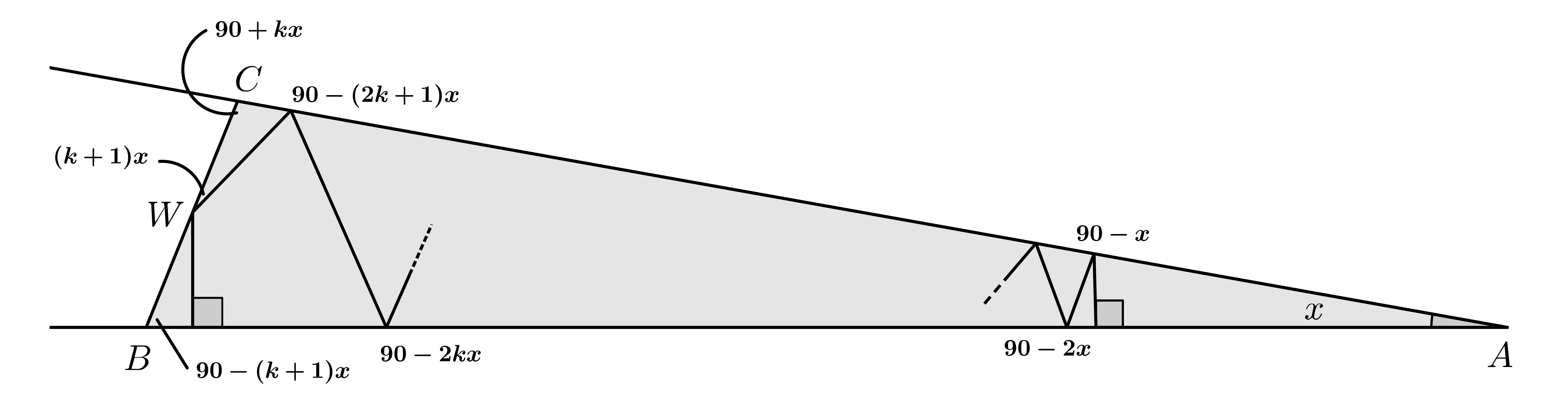}
        \caption{Infinite Pattern 2 Even Case}
        \label{fig:infd}
    \end{subfigure}
    ~
    \begin{subfigure}{1\textwidth}
        \centering
        \includegraphics[scale=0.18]{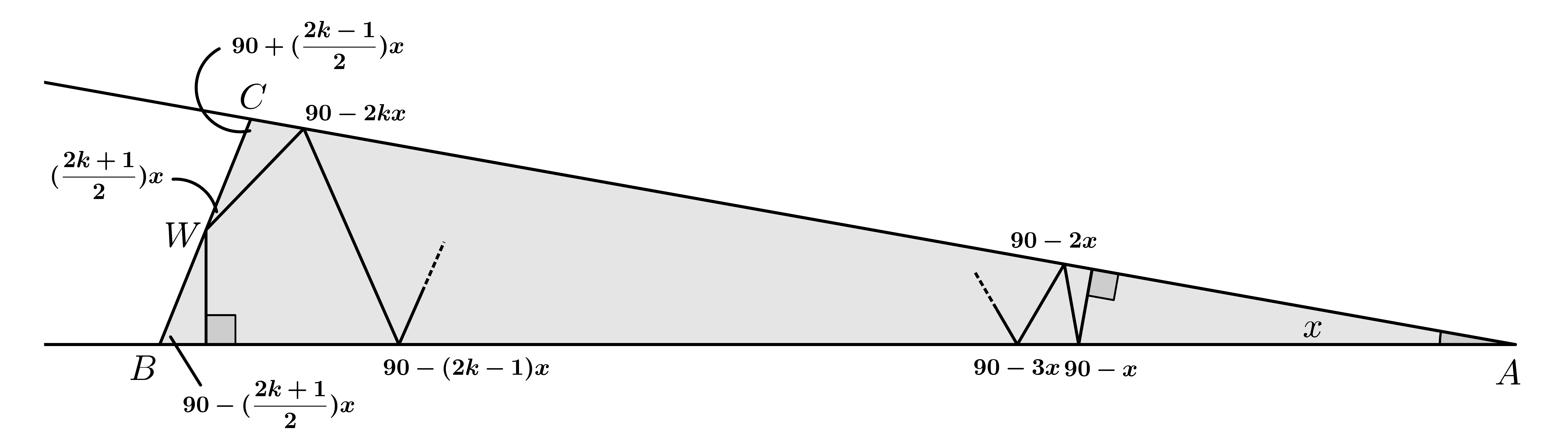}
        \caption{Infinite Pattern 2 Odd Case}
        \label{fig:infe}
    \end{subfigure}
    \caption{Infinite Pattern II}

\end{figure}

\begin{proof}
For the odd integer case $n=2k+1$, $k\geq0$, take a wedge of angle x and vertex A with arms $l_1$ and $l_2$ as shown in figure \ref{fig:infd} and shoot a poolball at 90 degrees from $l_1$ which then bounces off the sides at angles $90-x>90-2x>...>90-(2k+1)x>0$ noting that $(2k+1)x<90$. Now on the last ray leaving $l_2$ at angle $90-(2k+1)x$, choose any point W between $l_1$ and $l_2$ and draw a line through W hitting $l_1$ at B at an acute angle $y=90-(k+1)x>0$ (and so $(n+1)x+2y=180$) and hitting $l_2$ at C. Observe that triangle ABC has a periodic path of type 1 2 1 2$n$.

For the even integer case n=2k, $k\geq1$, again take a wedge of angle x and vertex A with arms $l_1$ and $l_2$ as shown in figure \ref{fig:infe} and shoot a poolball at 90 degrees from $l_2$ which then bounces off the sides at angles $90-x>90-2x>...>90-2kx>0$ noting that $2kx<90$. On the last ray leaving $l_2$ at angle $90-2kx$, choose any point W between $l_1$ and $l_2$ and draw a line through W hitting $l_1$ at B at an acute angle $y=90-(2k+1)x/2>0$ (and so $(n+1)x+2y=180$) and hitting $l_2$ at C. Observe that triangle ABC has a periodic path of type 1 2 1 2$n$.
\end{proof}

\begin{figure}[ht]
    \centering
	\begin{subfigure}[t]{0.5\textwidth}
    	\centering
		\includegraphics[scale=0.19]{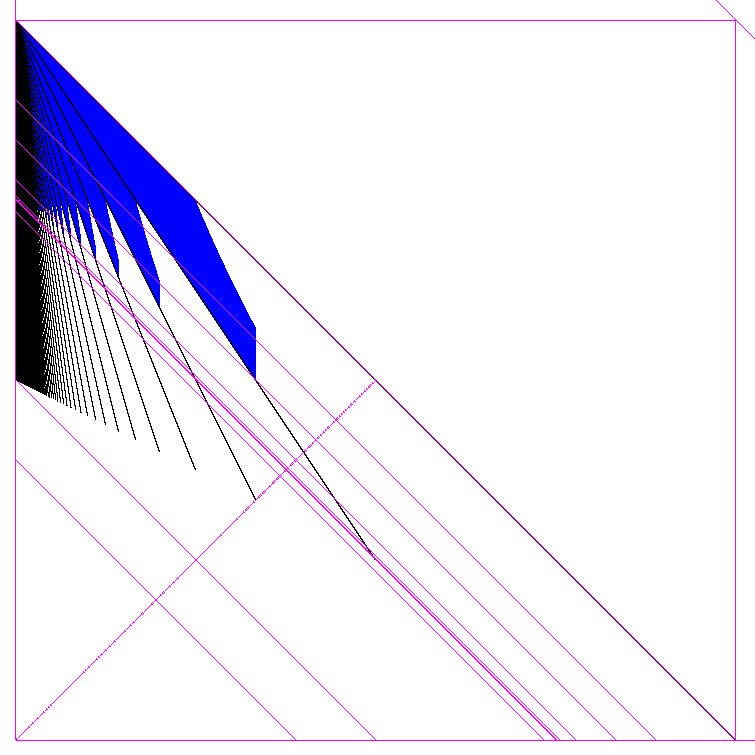}
		\caption{Infinite Patterns}
    \end{subfigure}%
    \begin{subfigure}[t]{0.5\textwidth}
    	\centering
    	\includegraphics[scale=0.12]{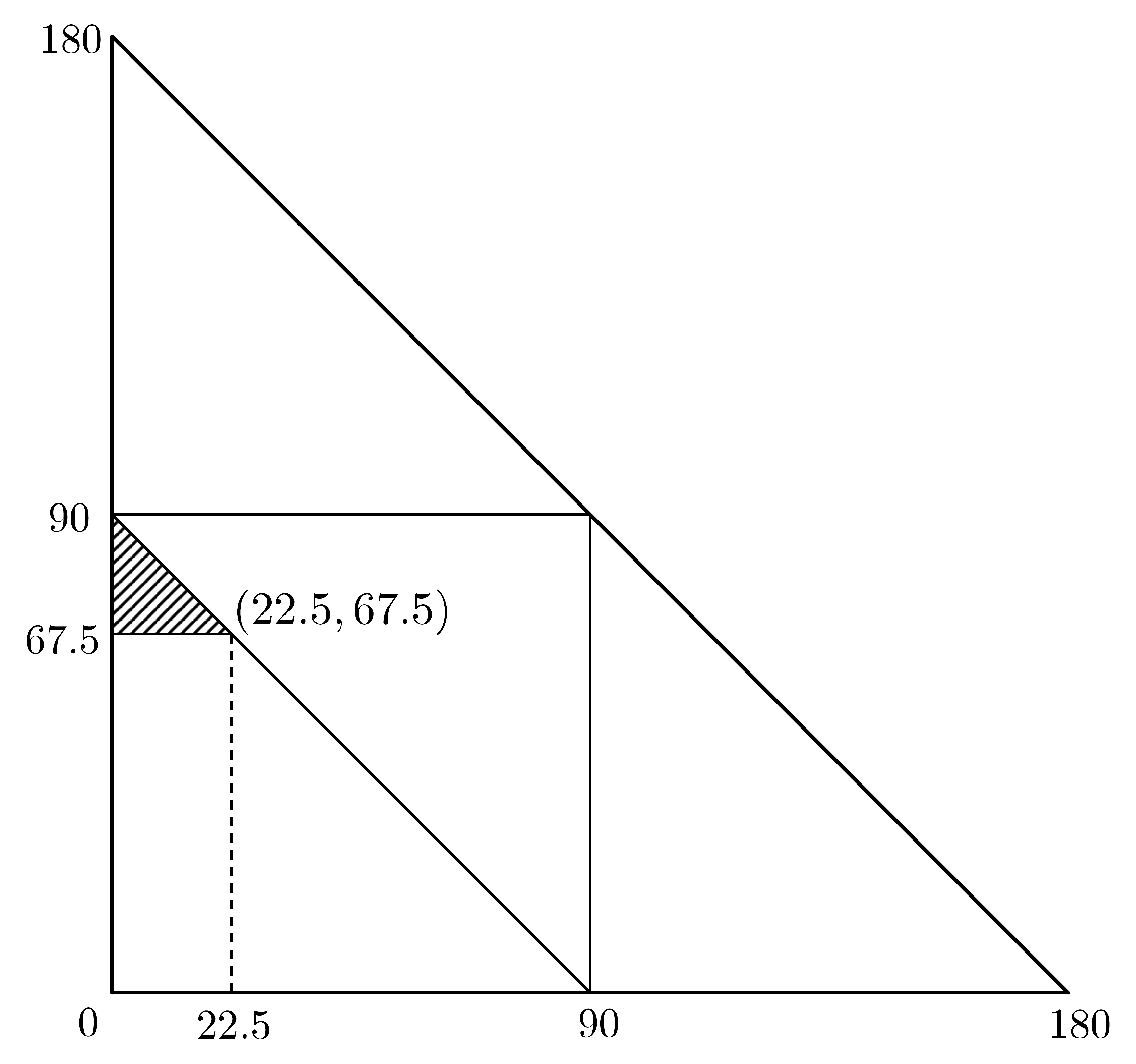}
    	\caption{Region Covered}
    \end{subfigure}
    \caption{Infinite Pattern Cover}
\end{figure}

\noindent
$\textbf{IMPORTANT CONCLUSION:}$ Given an obtuse triangle ABC, with $0<x<22.5$, $67.5<y<90$ and $x+y<90$, then that triangle has a periodic path.

\section{Bounding Polygons}

\noindent A \textbf{bounding polygon} is a convex polygon which includes a code region. Every code region has a bounding polygon for example the region bounded by $0<x+y<180$. It is useful in our calculations to find a bounding polygon with rational vertices which is as small as possible, the smallest being the convex hull of the region. There are up to six bounding polygons for each code, one corresponding to each permutation of the angles. For each code, we will assume we have a fixed order of the code angles $X$,$Y$ and $Z$.

\noindent
\textbf{The corner bounding polygon:} This is the polygon determined by the conditions that $0<nX<180$, $0<mY<180$ and $0<pZ<180$ where $nX,mY,pZ$ are the code angles corresponding to the largest $X$,$Y$,$Z$ code numbers in the code sequence.

\noindent \textbf{The angle bounding polygon:} Given a periodic side sequence leaving side AB of triangle ABC at an angle $T$ where $0<T \leq 90$, then we can calculate its successive angles as it reflects off each side. Since $z=180-x-y$, these reflecting angles will be linear combinations of $x$, $y$, 90 and $T$ with integer coefficients. If $T$ can be expressed in terms of $x$, $y$ and 90 with integer coefficients, then so can all reflecting angles. This is the case for the OSO, CS and CNS periodic paths but not for the OSNO or ONS periodic paths. An example is given below.

\begin{figure}[ht]
    \centering
    \begin{subfigure}[t]{0.33\textwidth}
        \centering
    	\includegraphics[scale=0.2]{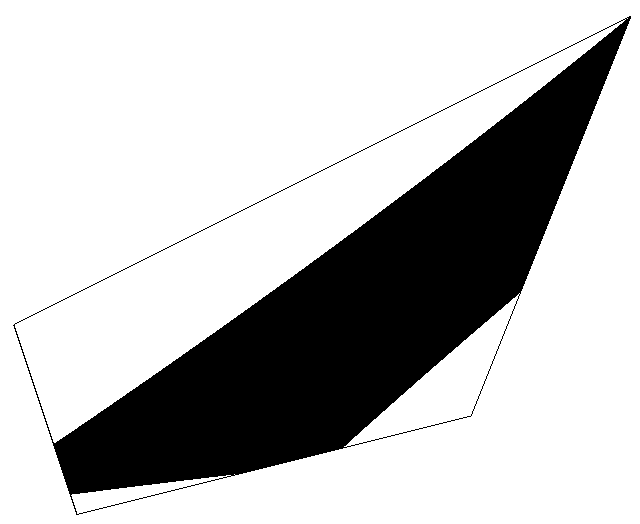}
    	\caption{CS Code 1 3 6 3 1 7 3 1 8 1 3 7}
    \end{subfigure}%
    \begin{subfigure}[t]{0.33\textwidth}
        \centering
    	\includegraphics[scale=0.15]{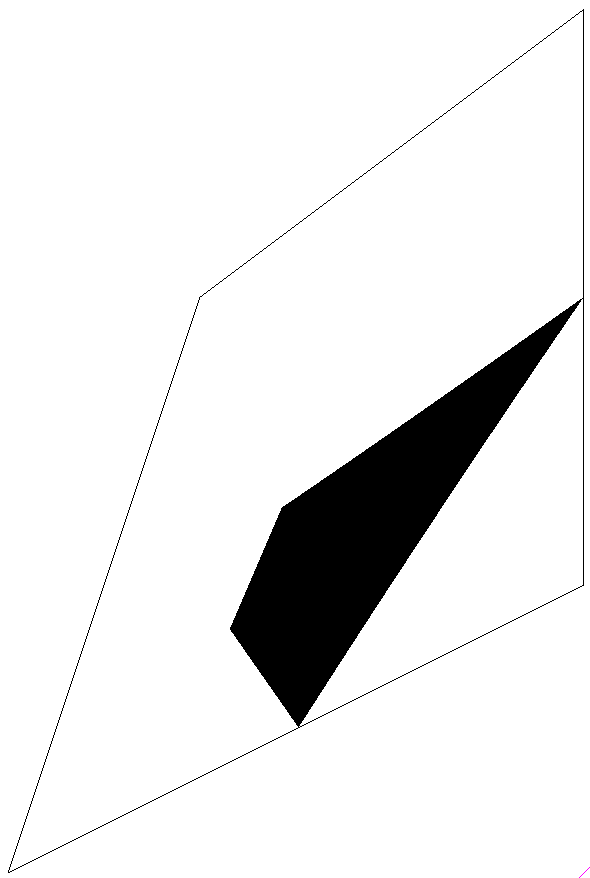}
    	\caption{OSO Code 1 1 2 2 3 1 2 1 4}
    \end{subfigure}
    \begin{subfigure}[t]{0.33\textwidth}
        \centering
    	\includegraphics[scale=0.2]{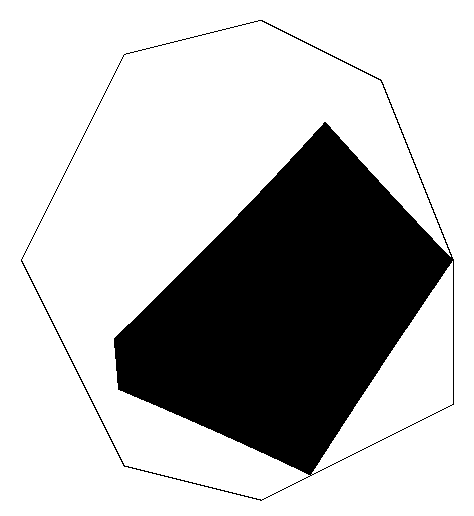}
    	\caption{OSNO Code 1 1 2 2 3 1 3 5}
    \end{subfigure}
    \caption{Angle Bounding Polygons}
\end{figure}

\noindent Now observe that in the OSO, CS and CNS cases each reflecting angle $\theta$ must satisfy $0<\theta\leq90$ and if we omit the 90 degree angles, then the set of $(x,y)$ satisfying $0<\theta<90$ forms a bounding polygon which contains the region determined by the given periodic path. In the CNS case, it is a bounding line segment.

\noindent On the other hand in the OSNO and ONS cases, we must have $0<\theta<90$ since there are no 90 degree angles. However since these reflecting angles are expressed in terms of $x$, $y$, 90 and $T$ with integer coefficients, in order to form the bounding polygon we must eliminate $T$. This can be done as follows. For each reflecting angle of the form $mx+ny+p90+T$, we can also say that $0<90-mx-ny-p90-T<90$ and similarly for reflecting angles of the form $mx+ny+p90-T$. We then get two sets of angles involving either $T$ or $-T$ and if we add each equation with $T$ to each equation with $-T$ and divide by 2, we end up with a set of linear cominations of $x$,$y$ and 90 with rational coefficients which lie between 0 and 90 and hence produce a bounding polygon in these cases. In the ONS case it is a bounding line segment. These are the bounding polygons
that we usually use in our calculations.

It is worth noting that the corner bounding polygon equations are included amongst the angle bounding polygon equations. This is a consequence of the fact that if a poolshot enters a corner $A$ where $m<A=x$ at an angle $\theta$ and bounces n times before it leaves then the angles involved are $\theta$, $\theta +x$, $\theta +2x$, ... , $180- \theta - (n-2)x$, $180- \theta - (n-1)x$ , $180- \theta - nx$ and then since $0< \theta <90$ and $0< 180- \theta - nx <90$, we must have $0< 180- nx <180$ or $0< nx <180$ which is one of the corner equations.

\section{The Program and Proof}

Because of the complexity and quantity of these equations, there is a dire need to automate the process of proving the codes work. Thus, we have written a program to crunch the numbers. In these calculations, each $G_j$ is an integer and is exact. On the other hand $r$, $b_{1}$, and $b_{2}$ are in radians and in fact are rational multiples of $\pi/2$, so they too are exact. The $f_j$ involve evaluating sines and cosines, so they are not exact. All these calculations are done by computer and have a certain degree of accuracy. We need to make sure that when we calculate that $f_j>0$, that it is indeed true. To do this we use interval arithmetic which can show that $f_j$ lies exactly within an interval $[u,v]$ with $u>0$. The interval that we use for $\pi/2$ correct to 7 decimal places is (1.57079631,1.57079637). This precision can be increased as required.

\noindent
We mainly use the arithmetical operations of\newline

addition: $[x_1,x_2]+[y_1,y_2]=[x_1+y_1,x_2+y_2]$

subtraction: $[x_1,x_2]-[y_1,y_2]=[x_1-y_2,x_2-y_1]$

multiplication: $[x_1,x_2][y_1,y_2]=[min(x_1y_1,x_1y_2,x_2y_1,x_2y_2),max(x_1y_1,x_1y_2,x_2y_1,x_2y_2)]$\newline

\noindent Note, we don't use any division operations.\newline

\noindent
In doing the calculations in the prover, we rely on the following libraries for arbitrary precision arithmetic.

\noindent
1. GMP: https://gmplib.org/

\noindent
2. MPFR: http://www.mpfr.org/ \cite{mpfr}

\noindent
We also rely on the boost multiprecision interface for all three libraries.\newline

\noindent
Along with this paper, you should get our program and instructions here using the terminal: \href{https://bitbucket.org/gtokarsky/billiardviewersm/src/master/}{Billiard Viewer}\newline

In the instructions of this program, we detail how you can view the proof. You can go square by square to see all the equations that produce the given code region and the lower bound using interval arithmetic of the calculations that show a square satisfies the prover. Because of how massive the proof is, there isn't a good way to present all of the calculations at once, so to trust this program is doing the calculations properly, we recommend that you read the code. 

In Appendix B, we have an ennumeration of the 134 code regions that cover the rest of the total region between $z=75$ and $z=80$ and which are not part of the two infinite patterns. The endpoints of this region in clockwise order are (37.5, 37.5),(40, 40),
(12.5, 67.5),(7.5, 67.5). We used 13,862 squares at 7 decimal accuracy starting from a single square and had to subdivide any square at most 20 times to be able to prove that these 134 code regions do cover this region. Our interval arithmetic calculations showed that the smallest lower bound on the prover at any equation used on a covering square was 6.65023 x $10^{-9}$.
The enumeration shows the code type, a 2-tuple consisting of the code length and the side sequence length followed by the code sequence. 

On this same program, you can find listed the 2439 single codes, 278,131 squares and 21 triples that are used to prove from 105 to 110, the 38,132 single codes,  4,994,538 squares and 310 triples to prove from 110 to 112 and the 118,809 single codes, 27,783,085 squares and 1,115 triples to prove from 112 to 112.3

\begin{figure}
	\centering
    \includegraphics[scale=0.3]{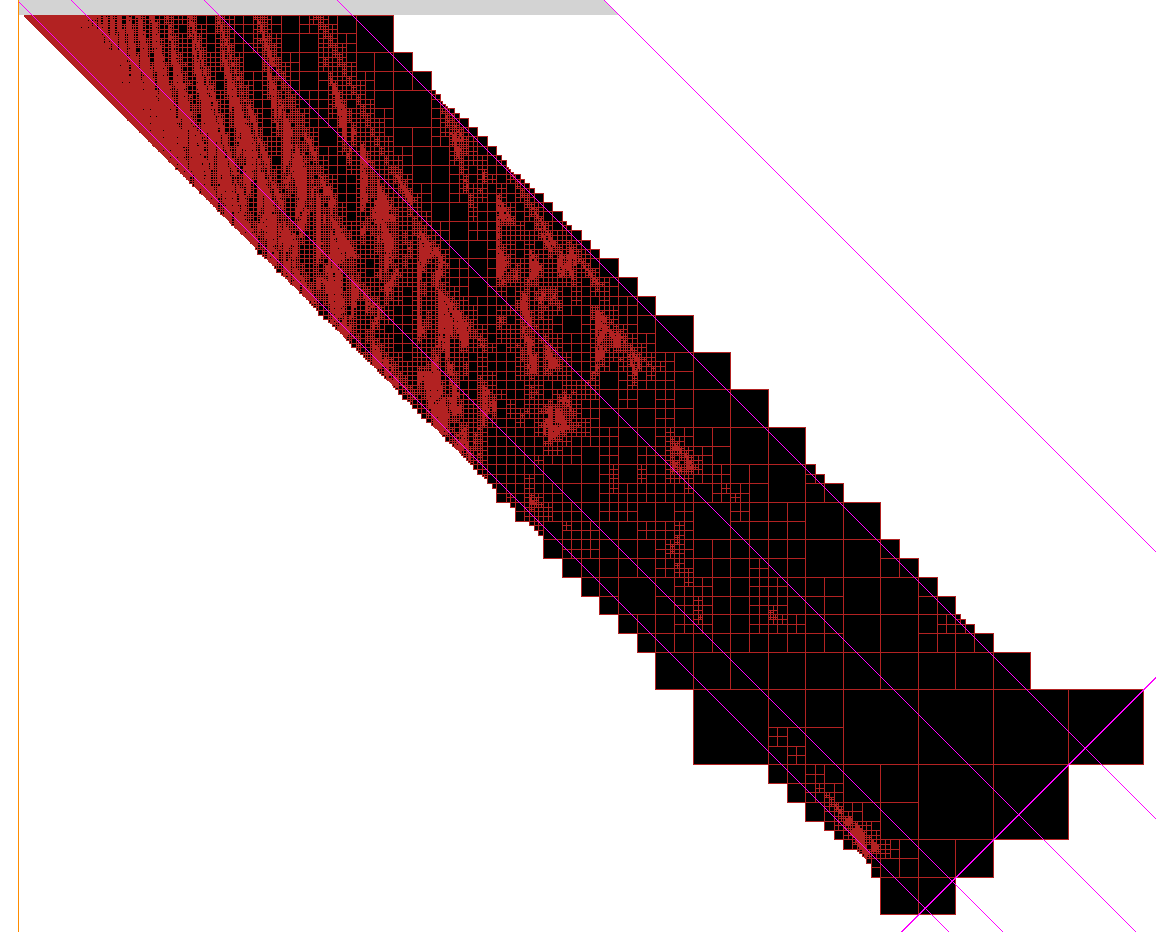}
    \caption{112.3 Cover}
\end{figure}

\section*{Acknowledgements}

The authors would like to thank the NSERC USRA program, Gerald Cliff, Brendan Pass, Eric Woolgar, and Terry Gannon, all of whom provided financial support to this project.

\begin{appendices}

	\section{Example of Calculating Coordinates of a Code Tower}
	
Here we will calculate the vertices in the code tower of the code ONS code sequence along the axis $x=y$.
\newline
\newline
	X~~~Z~~~X
	
\noindent
	1 1 2 3 3 2

\noindent
	~~Y~~~Y~~~Z

\noindent
Then

\noindent
$x_{1}=siny$, $y_{1}=0$ which are the coordinates of $L_{(1,0)}$

\noindent
$x_{2}$=0, $y_{2}=0$ which are the coordinates of $L_{(2,0)}$

\noindent
$x_{3}=sinzcosX$,

\noindent
$y_{3}=sinzsinX$ which are the coordinates of $L_{(3,0)}$

\noindent
$x_{4}=sinzcosX-sinxcos(Y-X)$,

\noindent$y_{4}=sinzsinX+sinxsin(Y-X)$ which are the coordinates of $L_{(4,0)}$

\noindent
$x_{5}=sinzcosX-sinxcos(Y-X)+sinxcos(2Z+X-Y)$,

\noindent$y_{5}=sinzsinX+sinxsin(Y-X)+sinxsin(2Z+X-Y)$ which are the coordinates of $L_{(5,0)}$

\noindent
$x_{6}=sinzcosX-sinxcos(Y-X)+sinxcos(2Z+X-Y)-sinzcos(3Y-2Z+Y-X)$,

\noindent$y_{6}=sinzsinX+sinxsin(Y-X)+sinxsin(2Z+X-Y)+sinzsin(3Y-2Z+Y-X)$ which are the coordinates of $L_{(6,0)}$

\noindent
$x_{7}=sinzcosX-sinxcos(Y-X)+sinxcos(2Z+X-Y)-sinzcos(3Y-2Z+Y-X)+sinycos(3X-3Y+2Z-Y+X)$,

\noindent$y_{7}=sinzsinX+sinxsin(Y-X)+sinxsin(2Z+X-Y)+sinzsin(3Y-2Z+Y-X)+sinysin(3X-3Y+2Z-Y+X)$ which are the coordinates of $L_{(7,0)}$

\noindent
$x_{8}=sinzcosX-sinxcos(Y-X)+sinxcos(2Z+X-Y)-sinzcos(3Y-2Z+Y-X)+sinycos(3X-3Y+2Z-Y+X)-sinycos(2Z-3X+3Y-2Z+Y-X)$,

\noindent$y_{8}=sinzsinX+sinxsin(Y-X)+sinxsin(2Z+X-Y)+sinzsin(3Y-2Z+Y-X)+sinysin(3X-3Y+2Z-Y+X)+sinysin(2Z-3X+3Y-2Z+Y-X)$ which are the coordinates of $L_{(8,0)}$
\newline
\newline
\noindent We then use standard trig identities together with the conditions that $z=180-x-y$ and $y=x$ from the code pattern.
\newline

$sin(-A)=-sinA$

$cos(-A)=cosA$

$sin(180-A)=sinA$

$cos(180-A)=-cosA$

$sinA+sinB=2sin(\frac{A+B}{2})cos(\frac{A-B}{2})$

$sinA-sinB=2sin(\frac{A-B}{2})cos(\frac{A+B}{2})$

$cosA+cosB=2cos(\frac{A+B}{2})cos(\frac{A-B}{2})$

$cosA-cosB=-2sin(\frac{A+B}{2})sin(\frac{A-B}{2})$

$2cosAcosB=cos(A+B)+cos(A-B)$

$2sinAcosB=sin(A+B)+sin(A-B)=sin(A+B)-sin(B-A)$

$2sinAsinB=cos(A-B)-cos(A+B)$
\newline\newline
\noindent
The shooting vector is $(c,d)$ where $(c,d)$ is the vector from $L_{(2,0)}$ to $L_{(8,0)}$

\noindent
$c= sin(y)cos(180)+sin(z)cos(x)+sin(x)cos(x-y+180)+sin(x)cos(-x-3y+4(180))+sin(z)cos(-x-6y+5(180))+sin(y)cos(2x-6y+6(180))$

\noindent
$d= sin(y)sin(180)+sin(z)sin(x)+sin(x)sin(x-y+180)+sin(x)sin(-x-3y+4(180))+sin(z)sin(-x-6y+5(180))+sin(y)sin(2x-6y+6(180))$

\noindent
\textbf{CONVENTION:} In using the last three trig identities to simplify to sums of sines and cosines, we multiply all coordinates by 2 so that all coefficients are integers.

\noindent
The shooting vector $(c,d)$ then further becomes

\noindent
$c= -2sin(y)-sin(3y)+sin(5y)-sin(2x-7y)+sin(2x-5y)-sin(2x-y)+sin(2x+y)+sin(2x+3y)-sin(2x+7y)$ a sum of sines.

\noindent
$d=-cos(3y)+cos(5y)+cos(2x-7y)-cos(2x-5y)+cos(2x-y)-cos(2x+y)+cos(2x+3y)-cos(2x+7y)$ a sum of cosines.

\noindent To calculate the coordinates of the key points in the tower which are not centers of fans if any, we illustrate by the same example. Suppose we look at the coordinates of $L_{(4,1)}$. It is found by starting with the coordinates of $L_{(4,0)}$ and since the corresponding code is 2 adding one more reflection of the given triangle and proceeding as before.

\noindent
$x_{4}=sinzcosX-sinxcos(Y-X)$,

\noindent$y_{4}=sinzsinX+sinxsin(Y-X)$

becomes

\noindent
$x_{(4,1)}=sinzcosX-sinxcos(Y-X)+ sinycos(Z-Y+X)$

\noindent
$y_{(4,1)}=sinzsinX+sinxsin(Y-X)+ sinysin(Z-Y+X)$

\noindent and then simplyfing to sums of sines or cosines.

\noindent To find the coordinates of $L_{(6,2)}$, we would start with the coordinates of $L_{(6,0)}$ and since the corresponding code is 3 adding two more reflections of the given triangle which corresponds to using the angle $2X-3Y+2Z-Y+X$ and we get

\noindent
$x_{(6,2)}=x_{6}+sinzcos(2X-3Y+2Z-Y+X )=x_{6}+sinzcos(x-6y)$

\noindent
$y_{(6,2)}=y_{6}+sinzsin(2X-3Y+2Z-Y+X )=x_{6}+sinzsin(x-6y)$

\noindent and then simplyfing to sums of sines or cosines.

\noindent
It is now a simple matter to calculate a vector $(a,b)$ from any key blue point to any key black point. As an example the blue-black vector from $L_{(6,0)}$ to $L_{(5,0)}$ is given by

$a=sinzcos(3Y-2Z+Y-X)=sinzcos(6y+x)=sin(7y+2x)-sin5x=sin9x-sin5x$ since y=x for this code

$b=-sinzsin(3Y-2Z+Y-X)=-sinzsin(6y+x)=cos(7y+2x)-cos5y=cos9x-cos5x$

\section{The Single Codes for the 105 theorem}

\begin{enumerate}[noitemsep]
\item OSO (3, 7) 1 3 3
\item OSO (5, 11) 1 1 2 2 5
\item OSO (5, 15) 1 1 4 2 7
\item OSO (5, 15) 1 3 2 6 3
\item OSO (5, 17) 1 1 4 2 9
\item OSO (5, 21) 1 1 6 2 11
\item OSO (5, 23) 1 1 6 2 13
\item OSO (7, 15) 1 1 3 1 2 1 6
\item OSO (7, 17) 1 1 3 1 2 1 8
\item OSO (7, 19) 1 1 2 2 6 2 5
\item OSO (7, 21) 1 1 2 2 8 2 5
\item OSO (7, 23) 1 1 4 2 6 2 7
\item OSO (7, 29) 1 1 6 2 8 2 9
\item OSO (7, 17) 1 2 1 2 1 3 7
\item OSNO (8, 18) 1 1 2 2 3 1 3 5
\item OSNO (8, 22) 1 1 4 2 3 1 3 7
\item OSO (9, 25) 1 1 2 2 6 2 4 2 5
\item CS (10, 20) 1 1 3 1 2 1 3 1 1 6
\item OSNO (10, 24) 1 1 3 1 2 1 5 1 1 8
\item CS (10, 28) 1 1 5 1 2 1 5 1 1 10
\item OSNO (10, 32) 1 1 5 1 2 1 7 1 1 12
\item CS (10, 36) 1 1 7 1 2 1 7 1 1 14
\item OSNO (10, 40) 1 1 7 1 2 1 9 1 1 16
\item CS (10, 44) 1 1 9 1 2 1 9 1 1 18
\item OSNO (10, 48) 1 1 9 1 2 1 11 1 1 20
\item OSNO (10, 26) 1 1 2 2 7 1 1 4 2 5
\item OSNO (10, 38) 1 1 4 2 11 1 1 6 2 9
\item OSO (11, 29) 1 1 2 1 1 7 2 4 1 1 8
\item OSO (11, 25) 1 1 2 1 1 6 1 1 2 2 7
\item OSNO (12, 38) 1 1 2 2 8 2 3 1 3 8 2 5
\item OSNO (14, 40) 1 1 3 1 2 1 7 2 4 1 1 7 2 7
\item OSNO (14, 44) 1 1 3 1 2 1 7 2 6 1 1 9 2 7
\item OSNO (14, 52) 1 1 3 1 2 1 11 2 6 1 1 11 2 9
\item OSNO (14, 36) 1 1 2 2 7 1 2 1 3 1 1 7 2 5
\item CS (14, 36) 1 2 1 5 3 1 3 2 4 2 3 1 3 5
\item CS (16, 44) 1 1 2 1 1 7 3 1 3 2 6 2 3 1 3 7
\item CS (16, 52) 1 1 4 1 1 9 3 1 3 2 8 2 3 1 3 9
\item OSNO (18, 42) 1 1 2 1 1 6 1 1 2 2 7 1 2 1 2 1 3 7
\item CS (18, 44) 1 1 2 2 5 1 2 1 5 2 2 1 1 5 2 4 2 5
\item CS (18, 52) 1 1 4 2 5 1 2 1 5 2 4 1 1 7 2 4 2 7
\item OSO (19, 57) 1 1 2 1 1 6 1 2 1 3 1 1 7 2 8 2 8 2 7
\item CS (20, 48) 1 1 3 1 2 1 6 1 2 1 3 1 1 8 1 1 4 1 1 8
\item CS (20, 56) 1 1 3 1 2 1 10 1 2 1 3 1 1 10 1 1 4 1 1 10
\item CS (20, 64) 1 1 5 1 2 1 8 1 2 1 5 1 1 12 1 1 6 1 1 12
\item CS (20, 88) 1 1 9 1 2 1 8 1 2 1 9 1 1 18 1 1 10 1 1 18
\item OSNO (20, 74) 1 1 2 1 1 7 1 1 13 2 7 1 2 1 7 1 1 13 2 9
\item CS (20, 52) 1 1 2 1 1 7 2 2 1 1 5 2 6 2 5 1 1 2 2 7
\item CS (20, 60) 1 1 2 1 1 7 2 4 1 1 7 2 6 2 7 1 1 4 2 7
\item CS (20, 68) 1 1 4 1 1 9 2 4 1 1 7 2 8 2 7 1 1 4 2 9
\item CS (20, 92) 1 1 6 1 1 13 2 6 1 1 11 2 10 2 11 1 1 6 2 13
\item CS (20, 60) 1 1 2 2 9 1 1 4 2 8 2 4 1 1 9 2 2 1 1 6
\item CS (20, 92) 1 1 6 2 13 1 1 8 2 10 2 8 1 1 13 2 6 1 1 12
\item CS (20, 60) 1 1 4 2 6 2 4 1 1 7 2 4 1 1 8 1 1 4 2 7
\item OSO (21, 49) 1 1 1 1 2 1 7 2 5 1 1 2 1 1 7 2 5 1 2 1 4
\item OSNO (22, 54) 1 1 1 1 3 7 1 1 3 1 2 1 5 2 6 2 2 1 1 5 2 5
\item OSNO (22, 82) 1 1 2 1 1 7 1 1 12 1 1 4 1 1 11 2 8 1 1 13 2 9
\item OSNO (22, 58) 1 1 2 2 7 1 2 1 3 2 7 1 1 3 1 2 1 5 2 6 2 5
\item OSNO (22, 64) 1 1 2 2 7 1 2 1 3 2 6 2 4 1 1 7 2 4 2 6 2 5
\item OSNO (22, 60) 1 1 4 2 7 1 2 1 4 1 2 1 5 1 1 9 2 3 1 2 1 8
\item OSNO (22, 74) 1 1 2 2 8 2 4 2 6 2 7 1 1 4 2 6 2 6 2 6 2 5
\item OSO (23, 55) 1 1 2 1 1 5 1 1 8 1 2 1 2 1 2 1 8 1 1 4 1 1 8
\item CS (24, 64) 1 1 1 1 2 1 9 3 1 2 1 3 9 1 2 1 1 1 1 5 2 8 2 5
\item OSNO (24, 60) 1 1 3 1 2 1 5 2 3 1 3 6 2 3 1 3 7 1 1 3 1 2 1 6
\item OSNO (24, 88) 1 1 3 1 2 1 11 2 9 1 2 1 5 1 1 12 1 1 6 1 1 13 2 9
\item OSNO (24, 68) 1 1 3 1 2 1 8 1 1 4 1 1 9 2 4 1 1 9 2 3 1 2 1 8
\item CS (24, 64) 1 1 2 1 1 5 1 1 10 1 1 6 1 1 10 1 1 5 1 1 2 1 1 8
\item CS (24, 76) 1 1 4 2 6 2 3 1 3 8 3 1 3 2 6 2 4 1 1 7 2 4 2 7
\item CS (26, 64) 1 1 1 1 2 1 7 2 4 2 7 1 2 1 1 1 1 5 2 5 1 2 1 5 2 5
\item CS (26, 104) 1 1 5 1 2 1 9 2 12 2 9 1 2 1 5 1 1 13 2 7 1 2 1 7 2 13
\item OSNO (26, 72) 1 1 2 1 1 5 1 1 8 1 2 1 3 2 9 1 1 4 2 9 1 1 4 1 1 8
\item CS (26, 64) 1 1 2 1 1 6 1 1 2 1 1 7 3 1 3 2 5 1 2 1 5 2 3 1 3 7
\item CS (26, 124) 1 1 6 2 12 2 7 1 2 1 7 2 12 2 6 1 1 11 2 8 2 12 2 8 2 11
\item CS (28, 80) 1 1 3 1 2 1 7 2 4 1 1 8 1 1 4 2 7 1 2 1 3 1 1 7 2 6 2 7
\item CS (28, 112) 1 1 3 1 2 1 11 2 8 1 1 14 1 1 8 2 11 1 2 1 3 1 1 9 2 12 2 9
\item OSNO (28, 118) 1 1 3 1 2 1 11 2 8 2 10 2 11 1 1 6 2 12 2 8 2 11 1 1 4 1 1 10
\item CS (28, 76) 1 1 3 1 2 1 10 1 1 4 1 1 10 1 2 1 3 1 1 8 1 2 1 6 1 2 1 8
\item CS (28, 80) 1 1 2 1 1 5 1 1 9 2 4 2 9 1 1 5 1 1 2 1 1 9 3 1 2 1 3 9
\item OSNO (28, 74) 1 1 2 1 1 5 2 6 2 5 1 1 2 2 7 1 2 1 3 2 6 2 5 1 1 2 2 7
\item OSNO (28, 68) 1 1 2 2 6 2 5 1 1 2 2 7 1 1 3 1 2 1 6 1 1 3 1 2 1 5 2 5
\item CS (28, 88) 1 1 4 2 10 2 4 1 1 7 2 9 1 1 5 1 2 1 4 1 2 1 5 1 1 9 2 7
\item CS (28, 92) 1 1 4 2 6 2 6 2 6 2 4 1 1 7 2 6 2 4 1 1 8 1 1 4 2 6 2 7
\item OSNO (30, 76) 1 1 2 2 5 1 2 1 4 1 2 1 5 3 1 2 1 3 6 3 1 3 2 4 2 4 2 4 2 5
\item OSO (31, 87) 1 1 3 1 2 1 6 1 2 1 3 1 1 7 2 6 2 4 1 1 7 2 5 1 2 1 5 2 6 2 7
\item OSNO (32, 84) 1 1 2 1 1 5 1 1 8 1 1 2 1 1 5 2 6 2 7 1 1 4 2 7 1 1 2 1 1 5 2 7
\item CS (32, 76) 1 1 2 1 1 5 1 1 8 1 2 1 2 1 2 1 8 1 1 5 1 1 2 1 1 8 1 1 4 1 1 8
\item CS (32, 76) 1 1 2 1 1 6 1 2 1 3 1 1 8 1 1 4 1 1 8 1 1 4 1 1 8 1 1 3 1 2 1 6
\item CS (32, 80) 1 1 2 2 5 1 2 1 5 3 1 2 1 3 5 1 2 1 5 2 2 1 1 5 2 4 2 4 2 4 2 5
\item OSNO (34, 80) 1 1 1 1 3 7 1 1 3 1 2 1 5 2 5 1 2 1 3 2 7 1 1 3 1 2 1 5 2 5 1 2 1 4
\item CS (34, 96) 1 1 3 1 2 1 5 2 6 2 4 2 6 2 5 1 2 1 3 1 1 7 2 4 2 5 1 2 1 5 2 4 2 7
\item CS (36, 84) 1 1 2 1 1 5 1 1 8 1 1 2 1 1 6 1 1 2 1 1 8 1 1 5 1 1 2 1 1 8 1 1 4 1 1 8
\item CS (36, 108) 1 1 2 1 1 7 1 1 12 1 1 4 1 1 10 1 1 4 1 1 12 1 1 7 1 1 2 1 1 10 1 1 4 1 1 10
\item CS (36, 132) 1 1 4 1 1 9 1 1 14 1 1 6 1 1 12 1 1 6 1 1 14 1 1 9 1 1 4 1 1 12 1 1 6 1 1 12
\item CS (36, 84) 1 1 2 1 1 6 1 1 2 1 1 6 1 1 2 1 1 7 3 1 3 2 5 1 2 1 4 1 2 1 5 2 3 1 3 7
\item CS (36, 84) 1 1 2 1 1 6 1 1 2 2 7 1 2 1 2 1 3 7 1 1 2 1 1 7 3 1 2 1 2 1 7 2 2 1 1 6
\item OSNO (36, 102) 1 1 2 1 1 6 1 2 1 5 2 7 1 1 3 1 2 1 6 1 2 1 5 2 6 2 6 2 7 1 1 4 2 6 2 7
\item CS (36, 112) 1 1 2 2 9 1 1 4 1 1 9 3 1 3 2 8 2 4 2 8 2 3 1 3 9 1 1 4 1 1 9 2 2 1 1 6
\item CS (38, 100) 1 1 2 1 1 6 1 1 2 1 1 7 2 7 1 1 3 1 2 1 7 2 5 1 2 1 5 2 7 1 2 1 3 1 1 7 2 7
\item OSNO (40, 96) 1 1 1 1 3 7 1 1 3 1 2 1 5 2 5 1 1 2 2 7 1 1 3 1 2 1 5 2 6 2 2 1 1 5 2 5 1 2 1 4
\item OSNO (40, 128) 1 1 3 1 2 1 10 1 1 5 1 2 1 7 2 13 1 1 7 1 2 1 6 1 2 1 7 2 13 1 1 7 1 2 1 6 1 2 1 8
\item OSNO (40, 120) 1 1 3 1 2 1 8 1 1 4 2 8 2 5 1 2 1 5 1 1 9 2 5 1 2 1 5 2 8 2 4 1 1 9 2 3 1 2 1 8
\item CS (40, 148) 1 1 2 1 1 7 1 1 13 2 9 1 2 1 5 1 1 13 2 8 2 13 1 1 5 1 2 1 9 2 13 1 1 7 1 1 2 1 1 10
\item CS (40, 120) 1 1 2 1 1 7 2 4 1 1 7 2 6 2 6 2 6 2 7 1 1 4 2 7 1 1 2 1 1 7 2 4 1 1 8 1 1 4 2 7
\item CS (40, 100) 1 1 2 1 1 6 1 1 3 1 2 1 5 2 6 2 5 1 2 1 3 1 1 6 1 1 2 1 1 7 2 4 1 1 8 1 1 4 2 7
\item OSNO (42, 116) 1 1 2 1 1 5 1 1 9 2 4 2 9 1 1 4 2 9 1 1 4 1 1 8 1 2 1 3 1 1 8 1 1 2 1 1 6 1 1 2 2 9
\item OSNO (44, 118) 1 1 2 1 1 7 2 4 1 1 8 1 1 2 1 1 7 2 4 1 1 8 1 2 1 3 2 9 1 1 5 1 2 1 4 1 2 1 7 2 4 1 1 8
\item CS (48, 140) 1 1 1 1 3 9 1 1 2 1 1 5 1 1 9 2 4 2 9 1 1 4 1 1 9 2 4 2 9 1 1 5 1 1 2 1 1 9 3 1 1 1 1 5 2 8 2 5
\item CS (48, 136) 1 1 2 1 1 6 1 1 2 1 1 7 2 2 1 1 5 2 6 2 4 2 7 1 1 4 2 6 2 4 2 6 2 4 1 1 7 2 4 2 6 2 5 1 1 2 2 7
\item CS (50, 156) 1 1 1 1 3 9 1 1 4 2 8 2 5 1 2 1 5 2 8 2 5 1 2 1 5 2 8 2 4 1 1 9 3 1 1 1 1 5 2 8 2 5 1 2 1 5 2 8 2 5
\item OSNO (50, 156) 1 1 2 1 1 7 2 6 2 6 2 8 2 4 1 1 9 2 3 1 2 1 8 1 1 3 1 2 1 8 1 1 4 2 8 2 6 2 6 2 7 1 2 1 3 1 1 7 2 7
\item CS (56, 132) 1 1 1 1 2 1 7 2 4 1 1 8 1 1 4 2 7 1 2 1 1 1 1 4 1 2 1 5 2 7 1 1 3 1 2 1 6 1 1 2 1 1 6 1 2 1 3 1 1 7 2 5 1 2 1 4
\item CS (56, 160) 1 1 2 1 1 5 1 1 9 2 4 2 9 1 1 4 1 1 9 2 4 2 9 1 1 5 1 1 2 1 1 9 3 1 2 1 3 9 1 1 2 1 1 6 1 1 2 1 1 9 3 1 2 1 3 9
\item OSNO (56, 216) 1 1 4 1 1 9 2 11 1 1 5 1 2 1 8 1 2 1 5 1 1 12 1 1 6 1 1 13 2 6 2 13 1 1 6 2 13 1 1 6 1 1 12 1 1 4 1 1 9 2 12 2 6 1 1 12
\item CS (56, 216) 1 1 2 1 1 8 1 1 2 1 1 9 2 12 2 9 1 2 1 5 1 1 13 2 7 1 2 1 7 2 13 1 1 6 1 1 13 2 7 1 2 1 7 2 13 1 1 5 1 2 1 9 2 12 2 9
\item OSNO (56, 142) 1 1 2 1 1 6 1 2 1 3 1 1 7 2 5 1 2 1 4 1 2 1 6 1 1 2 1 1 7 2 5 1 1 2 2 8 2 5 1 1 2 2 9 1 1 5 1 2 1 4 1 2 1 5 2 7
\item OSNO (56, 170) 1 1 2 1 1 6 1 2 1 5 2 8 2 4 1 1 8 1 2 1 3 1 1 9 2 2 1 1 5 1 1 9 2 4 1 1 9 2 4 2 8 2 6 2 6 2 7 1 2 1 3 1 1 7 2 7
\item CS (56, 172) 1 1 2 2 9 1 1 4 1 1 9 2 2 1 1 6 1 1 2 2 9 1 1 4 1 1 9 3 1 3 2 8 2 4 2 8 2 4 2 8 2 3 1 3 9 1 1 4 1 1 9 2 2 1 1 6
\item OSNO (58, 158) 1 1 2 1 1 6 1 2 1 5 2 7 1 2 1 3 1 1 7 2 7 1 1 2 1 1 6 1 2 1 5 2 8 2 4 1 1 9 2 2 1 1 5 1 1 9 2 4 1 1 8 1 2 1 3 1 1 8
\item OSNO (60, 162) 1 1 2 1 1 5 1 1 9 2 4 2 9 1 1 4 2 9 1 1 4 1 1 9 2 3 1 2 1 8 1 2 1 2 1 3 9 1 1 2 1 1 6 1 1 2 1 1 8 1 1 2 1 1 6 1 1 2 2 9
\item CS (60, 144) 1 1 2 1 1 5 1 1 8 1 2 1 2 1 2 1 8 1 1 4 1 1 8 1 2 1 2 1 2 1 8 1 1 5 1 1 2 1 1 8 1 1 4 1 1 8 1 2 1 2 1 2 1 8 1 1 4 1 1 8
\item CS (60, 180) 1 1 2 1 1 7 2 4 1 1 7 2 6 2 6 2 6 2 6 2 6 2 7 1 1 4 2 7 1 1 2 1 1 7 2 4 1 1 8 1 1 4 2 7 1 1 2 1 1 7 2 4 1 1 8 1 1 4 2 7
\item OSNO (60, 186) 1 1 2 1 1 7 2 6 2 6 2 8 2 4 1 1 8 1 2 1 3 1 1 8 1 2 1 3 2 9 1 1 4 2 8 2 6 2 6 2 7 1 1 2 1 1 7 2 6 2 7 1 2 1 3 1 1 7 2 7
\item CS (64, 156) 1 1 2 1 1 6 1 1 2 1 1 8 1 1 3 1 2 1 7 2 6 2 7 1 2 1 3 1 1 8 1 1 2 1 1 6 1 1 2 1 1 9 2 2 1 1 6 1 1 2 1 1 8 1 1 2 1 1 6 1 1 2 2 9
\item CS (66, 200) 1 1 1 1 3 9 1 1 2 1 1 5 1 1 9 2 4 2 9 1 1 4 1 1 9 2 4 2 9 1 1 4 1 1 9 2 4 2 9 1 1 5 1 1 2 1 1 9 3 1 1 1 1 5 2 8 2 5 1 2 1 5 2 8 2 5
\item CS (66, 208) 1 1 1 1 3 9 1 1 4 2 8 2 5 1 2 1 5 2 8 2 5 1 2 1 5 2 8 2 5 1 2 1 5 2 8 2 4 1 1 9 3 1 1 1 1 5 2 8 2 5 1 2 1 5 2 8 2 5 1 2 1 5 2 8 2 5
\item OSNO (66, 230) 1 1 2 1 1 5 1 1 9 2 4 1 1 9 2 4 2 8 2 6 2 6 2 8 2 4 1 1 9 2 4 2 8 2 6 2 6 2 8 2 4 2 9 1 1 4 1 1 9 2 4 1 1 9 2 4 2 9 1 1 5 1 1 2 2 9
\item CS (68, 188) 1 1 2 1 1 6 1 1 2 1 1 8 1 2 1 3 2 9 1 1 4 2 8 2 5 1 2 1 6 1 1 2 1 1 8 1 1 4 2 8 2 4 1 1 8 1 1 2 1 1 6 1 2 1 5 2 8 2 4 1 1 9 2 3 1 2 1 8
\item CS (72, 216) 1 1 2 1 1 7 2 4 1 1 8 1 1 2 1 1 7 2 4 2 8 2 4 2 8 2 4 2 7 1 1 2 1 1 8 1 1 4 2 7 1 1 2 1 1 7 2 5 1 2 1 5 2 8 3 1 3 2 8 2 3 1 3 8 2 5 1 2 1 5 2 7
\item OSNO (74, 220) 1 1 2 1 1 6 1 2 1 5 2 8 2 4 2 9 1 1 4 2 8 2 6 2 6 2 7 1 1 2 1 1 7 2 7 1 1 3 1 2 1 7 2 6 2 7 1 1 3 1 2 1 7 2 5 1 2 1 5 1 1 9 2 4 1 1 8 1 2 1 3 1 1 8
\item CS (80, 324) 1 1 2 1 1 7 2 13 1 1 8 2 10 2 9 1 2 1 5 1 1 12 1 1 6 2 13 1 1 6 1 1 12 1 1 6 2 12 2 8 2 12 2 6 1 1 12 1 1 6 1 1 13 2 6 1 1 12 1 1 5 1 2 1 9 2 10 2 8 1 1 13 2 7 1 1 2 1 1 10
\item CS (84, 260) 1 1 1 1 3 9 1 1 2 1 1 5 1 1 9 2 4 2 9 1 1 4 1 1 9 2 4 2 9 1 1 4 1 1 9 2 4 2 9 1 1 4 1 1 9 2 4 2 9 1 1 5 1 1 2 1 1 9 3 1 1 1 1 5 2 8 2 5 1 2 1 5 2 8 2 5 1 2 1 5 2 8 2 5
\item CS (92, 252) 1 1 1 1 3 9 1 1 5 1 1 2 1 1 8 1 1 4 1 1 8 1 2 1 2 1 2 1 8 1 1 4 1 1 8 1 1 2 1 1 5 1 1 9 3 1 1 1 1 5 2 6 2 7 1 2 1 2 1 2 1 8 1 1 4 2 8 2 4 2 8 2 4 2 8 2 4 1 1 8 1 2 1 2 1 2 1 7 2 6 2 5
\item CS (100, 300) 1 1 2 1 1 7 2 4 1 1 7 2 6 2 6 2 6 2 6 2 6 2 6 2 6 2 6 2 6 2 7 1 1 4 2 7 1 1 2 1 1 7 2 4 1 1 8 1 1 4 2 7 1 1 2 1 1 7 2 4 1 1 8 1 1 4 2 7 1 1 2 1 1 7 2 4 1 1 8 1 1 4 2 7 1 1 2 1 1 7 2 4 1 1 8 1 1 4 2 7
\item CS (102, 320) 1 1 1 1 3 9 1 1 2 1 1 5 1 1 9 2 4 2 9 1 1 4 1 1 9 2 4 2 9 1 1 4 1 1 9 2 4 2 9 1 1 4 1 1 9 2 4 2 9 1 1 4 1 1 9 2 4 2 9 1 1 5 1 1 2 1 1 9 3 1 1 1 1 5 2 8 2 5 1 2 1 5 2 8 2 5 1 2 1 5 2 8 2 5 1 2 1 5 2 8 2 5
\item OSNO (116, 326) 1 1 2 1 1 6 1 2 1 5 2 7 1 2 1 4 1 2 1 6 1 2 1 5 2 7 1 1 2 1 1 8 1 1 4 2 8 2 5 1 2 1 5 2 8 2 5 1 2 1 5 2 9 1 1 5 1 2 1 5 2 8 2 5 1 2 1 5 2 8 2 4 1 1 8 1 1 2 1 1 7 2 5 1 2 1 6 1 2 1 4 1 2 1 7 2 5 1 2 1 5 1 1 9 2 4 1 1 8 1 2 1 3 1 1 8
\end{enumerate}

\end{appendices}

\end{document}